\documentclass{siamltex}

\usepackage{graphicx,amsmath,amsfonts,amssymb}

\usepackage{soul} 
\usepackage{epstopdf}
\usepackage{booktabs}
\usepackage{boxedminipage}
\usepackage{booktabs}
\usepackage{algorithm}
\usepackage[noend]{algpseudocode}
\usepackage[mathscr]{eucal}
\usepackage{bm}
\usepackage{algpseudocode}
\usepackage{comment}
\usepackage{psfrag}
\usepackage{xcolor}
\usepackage{enumitem}
\usepackage{cite} 
\usepackage{wrapfig} 

\makeatletter
\def\BState{\State\hskip-\ALG@thistlm}
\makeatother

\newcommand*{\de}{\mathop{}\!\mathrm{d}}

\def\ba{\begin{aligned}}

\def\be{\begin{equation}}
\def\bSigma{{\bm\Sigma}}

\def\bthe{{\vec\theta}}
\def\bThe{{\bm\Theta}}
\def\bthes{{\bthe_\sss}}
\def\bThes{{\bThe_\SSS}}
\def\bx{{\bf\xxx}}

\def\cF{{\mathcal\FFF}}

\def\cN{{\mathcal\NNN}}
\def\cS{{\mathcal \SSS}}

\def\cU{{\mathcal\UUU}}
\def\cV{{\mathcal \VVV}}
\def\cvk{{|\cV_\kkk|}}

\def\ddd{{d}}

\def\ea{\end{aligned}}
\def\e{{\bm e}}
\def\ee{\end{equation}}

\def\el {\nonumber }
\def\f{{\bm f}}
\def\FFF{{F}}
\def\g{{\bm g}}

\def\hhh{{h}}

\def\kkk{{k}}
\def\KKK{{K}}
\def\lll{{\ell}}
\def\LLL{{L}}
\def\LLLk{{\LLL_\kkk}}
\def\MMM{{M}}
\def\n{{\bm n}}
\def\nnn{{n}}
\def\NNN{{N}}

\def\phph{{\phi}}
\def\psps{{\psi}}
\def\ppsi{{\bm\psi}}
\def\Psps{{\bm\Psi}}
\def\R{{\mathbb R}}

\def\sss{{s}}
\def\SSS{{S}}
\def\SM{{\bf S}}
\def\Upup{{\bm\Upsilon}}
\def\u{{\bm U}}
\def\uuu{{u}}
\def\UUU{{U}}
\def\uuul{{\uuu_\LLL}}
\def\uuulk{{\uuu_\LLLk}}
\def\uuuh{{\uuu_\NNN}}
\def\uuuhs{{\uuu_{\NNN,\sss}}}
\def\vvv{{v}}
\def\VVV{{V}}
\def\vvvl{{\vvv_\LLL}}
\def\VVVl{{\VVV_\LLL}}
\def\VVVlk{{\VVV_\LLLk}}

\def\VVVh{{\VVV_\NNN}}
\def\xxx{{x}}
\def\z{{\bm z}}

\begin{document}

\title{
A Localized Reduced-Order Modeling Approach \\ for PDEs with Bifurcating Solutions\thanks{AA supported by the US Department of Energy Office of Science grant DE-SC0009324. MG supported by the US Air Force Office of Scientific Research grant FA9550-15-1-0001 and US Department of Energy Office of Science grant DE-SC0009324. MH and GR supported by European Union Funding for Research and Innovation through the European Research Council project H2020 ERC CoG 2015 AROMA-CFD project 681447, P.I. Prof. G. Rozza). AQ supported by US National Science Foundation grant DMS-1620384. }}

\author{
Martin Hess\thanks{SISSA, International School for Advanced Studies, Mathematics Area, mathLab, Trieste, Italy \texttt{mhess@sissa.it}}
\and
Alessandro Alla\thanks{PUC-Rio, Department of Mathematics, Rio de Janeiro, Brazil \texttt{alla@mat.puc-rio.br}}
\and
Annalisa Quaini\thanks{Department of Mathematics, University of Houston,     
Houston, USA \texttt{quaini@math.uh.edu}}
\and
Gianluigi Rozza\thanks{SISSA, International School for Advanced Studies, Mathematics Area, mathLab, Trieste, Italy \texttt{gianluigi.rozza@sissa.it}}
\and
Max Gunzburger\thanks{Florida State University, Department of Scientific Computing, Tallahassee, USA \texttt{mgunzburger@fsu.edu}}
}

\maketitle

\begin{abstract}
Reduced-order modeling (ROM) commonly refers to the construction, based on a few solutions (referred to as snapshots) of an expensive discretized partial differential equation (PDE), and the subsequent application of low-dimensional discretizations of partial differential equations (PDEs) that can be used to more efficiently treat problems in control and optimization, uncertainty quantification, and other settings that require multiple approximate PDE solutions. In this work, a ROM is developed and tested for the treatment of nonlinear PDEs whose solutions bifurcate as input parameter values change. In such cases, the parameter domain can be subdivided into subregions, each of which corresponds to a different branch of solutions. Popular ROM approaches such as proper orthogonal decomposition (POD), results in a global low-dimensional basis that does no respect not take advantage of the often large differences in the PDE solutions corresponding to different subregions. Instead, in the new method, the k-means algorithm is used to cluster snapshots so that within cluster snapshots are similar to each other and are dissimilar to those in other clusters. This is followed by the construction of local POD bases, one for each cluster. The method also can detect which cluster a new parameter point belongs to, after which the local basis corresponding to that cluster is used to determine a ROM approximation. Numerical experiments show the effectiveness of the method both for problems for which bifurcation cause continuous and discontinuous changes in the solution of the PDE.
\end{abstract}

\maketitle

\section{Introduction}\label{sec:intro}

We consider the problem of determining a function $\uuu\in\VVV$ such that 
\be\label{abs-pro}
    \cN(\uuu,\vvv;\bthe) = \cF(\vvv;\bthe)\quad\forall\, \vvv\in\VVV, \qquad \Longleftarrow \mbox{\em PDE model}
\ee
where $\bthe\in \bThe$ denotes a point in a parameter domain $\bThe\subset \R^\MMM$, $\VVV$ a function space, $\cN(\cdot,\cdot;\bthe)$ a given form that is linear in $\vvv$ but generally nonlinear in $\uuu$, and $\cF(\cdot)$ a linear functional on $\VVV$. Note that either or both $\cN$ and $\cF$ could depend on some or all the components of the parameter vector\footnote{For finite-dimensional spaces, we use the nomenclature ``points'' and ``vectors'' interchangeably.} $\bthe$. We view the problem \eqref{abs-pro} as a variational formulation of a nonlinear partial differential equation (PDE) or a system of such equations\footnote{In \S\ref{sec:nse}, we use the Navier-Stokes equations as a concrete setting to illustrate our methodology.} in which $\MMM$ parameters appear. However, it is not necessary to do so, e.g., we could treat strong forms of PDEs as well; the notation of \eqref{abs-pro} allows us to keep the exposition relatively simple.
 
We are interested in situations in which the solution $u$ of \eqref{abs-pro} differs in character for parameter vectors $\bthe$ in different subregions of the parameter domain $\bThe$. Such situations occur if $u$ undergoes bifurcations as $\bthe$ changes from one subregion to another; in this case we refer to the subdivision of $\bThe$ into subregions as a {\em bifurcation diagram}. We are particularly interested in situations that require solutions of \eqref{abs-pro} for a set of parameter vectors that span across two or more of the subregions of the bifurcation diagram. Such a situation arises in optimization problems for which an (iterative) optimization algorithm updates the parameter vector at each iteration so that it is possible that the updated parameter vector is in a different subregion than is the previous vector. Uncertainty quantification problems also give rise to the need to find solutions of \eqref{abs-pro} for possibly many parameter vectors that span across two or more subregions of the bifurcation diagram. 

In this paper, we assume that we have a priori knowledge of the bifurcation diagram. We do so because we focus on the construction and application of a reduced-order model (ROM) for problems with bifurcating solutions. This assumption only affects the first step in the construction of a ROM. In a follow-up paper we treat the case for which no knowledge of the bifurcation diagram is presumed. In \S\ref{sec:con}, we discuss what changes are needed so that our approach can handle the more general case.

In general, one cannot solve \eqref{abs-pro} for $u$ so that one instead chooses an approximating $\NNN$-dimensional subspace $\VVVh\subset\VVV$ and then seeks an approximation $\uuuh\in\VVVh$ satisfying the discretized system\footnote{In general, the forms $\cN$ and $\cF$ themselves are also discretized, e.g., because quadrature rules are used to approximate integrals appearing in their definition. However, here, we ignore such approximations, again to keep the exposition simple.}
\be\label{abs-dis}
    \cN(\uuuh,\vvv;\bthe) = \cF(\vvv;\bthe)\quad\forall\, \vvv\in\VVVh.\qquad \Longleftarrow\mbox{\em full-order model}
\ee
Note that if $\{\phph_\nnn(\bx)\}_{\nnn=1}^\NNN$ denotes a basis for $\VVVh$, we have that $\uuuh=\sum_{\nnn=1}^\NNN \uuu_\nnn \phph_\nnn(\bx)$ for some set of coefficients $\{\uuu_\nnn\}_{\nnn=1}^\NNN$.\footnote{Equation \eqref{abs-dis} represents a Galerkin type setting in which the trial function $\uuuh$ and test function $\vvv$ belong to the same space $\VVVh$. We could easily generalize our discussion to the Petrov-Galerkin case for which these function would belong to different approximating spaces.} 

Solving the discretized nonlinear PDE \eqref{abs-dis} for $\uuuh$ is often an expensive proposition, especially if multiple solutions are needed. For example, if $\VVVh$ denotes a finite element subspace constructed with respect to a meshing of nominal grid size $\hhh$ of a $\ddd$-dimensional spatial domain, we have that $\NNN=O(\hhh^{-\ddd})$ so that the discretized system \eqref{abs-dis} could be huge. For this reason, one is interested in building surrogates for the solution of \eqref{abs-dis} that are much less costly to evaluate so that obtaining approximations to the exact solution $\uuu$ of \eqref{abs-pro} for many choices of the parameter vector $\bthe\in\bThe$ now becomes feasible. Such surrogates are invariably constructed using a ``few'' solutions of the expensive, full-order discrete system \eqref{abs-dis} and can take on many forms such as interpolants and least-squares approximations. Here, we are interested in {\em reduced-order models} (ROMs) for which one constructs a relatively low-dimensional approximating subspace $\VVVl\subset\VVVh$ of dimension $\LLL$ that still contains an acceptably accurate approximation $\uuul$ to the solution $u$ of \eqref{abs-pro}. That approximation is determined from the reduced discrete system
\be\label{abs-romg}
    \cN(\uuul,\vvv;\bthe) = \cF(\vvv;\bthe)\quad\forall\, \vvvl\in\VVVl \qquad \Longleftarrow \mbox{\em global reduced-order model}
\ee
that, if $\LLL\ll\NNN$, is much cheaper to solve compared to \eqref{abs-dis}. We note that ROM systems are generally dense in the sense that, e.g., matrices associated with \eqref{abs-romg} are not sparse.

As already stated, the goal is to the determine ROM approximations $\uuul$ that are acceptably accurate approximations of the solutions $\uuu$ of the continuous model \eqref{abs-pro}. However, because the ROM approximation is constructed through the use of solutions $\uuuh$ of the spatially discretized problem \eqref{abs-dis}, we have two sources of error in the ROM approximation, an error that can be estimated (in a chosen norm $\|\cdot\|$) by
\be\label{abs-tri}
   \underbrace{\| \uuu - \uuul \|}_{
   \begin{tabular}{cc}
   \mbox{\em error in ROM}
   \\
   \mbox{\em approximation}
   \end{tabular}} 
   \le 
   \underbrace{\| \uuu - \uuuh \|}_{
   \begin{tabular}{cc}
   \mbox{\em error due to spatial}
   \\
   \mbox{\em approximation}
   \end{tabular}} 
 + \underbrace{\| \uuuh - \uuul \|.}_{
   \begin{tabular}{cc}
   \mbox{\em error due to ROM}
   \\
   \mbox{\em approximation of the}
   \\
   \mbox{\em spatial approximation}
   \end{tabular}}
\ee
We mostly focus on the last term in \eqref{abs-tri}. However, one should keep in mind that an efficient overall computational methodology should try to balance the two error terms in the right-hand side of \eqref{abs-tri} because there is not much sense in having the error due to the ROM be much smaller than the spatial error. This observation may be used to relate $\LLL$ to $\NNN$ so as to provide guidance as to what should be the dimension $\LLL$ of the ROM approximations space.

ROMs in the setting of bifurcating solutions are considered in the early papers \cite{NOOR1982955,Noor:1994,Noor:1983,NOOR198367} in the setting of buckling bifurcations in solid mechanics. More recently, in \cite{Terragni:2012} it is shown that a POD approach allows for considerable computational time savings for the analysis of bifurcations in some nonlinear dissipative systems. In \cite{Maday:RB2,PLA2015162}, a reduced basis (RB) method is used to track solution branches from bifurcation points arising in natural convection problems. A RB method is used in \cite{PR15} to investigate Hopf bifurcations in natural convection problems and in \cite{PITTON2017534} for symmetry breaking bifurcations in contraction-expansion channels.
An investigation of symmetry breaking in an expansion channel can be found in \cite{AQpreprint}. In \cite{Cliffe2010,Cliffe2011,cliffe}, reliable error estimation is used to determine the critical parameter points where bifurcations occur in the Navier-Stokes setting. A recent work on ROMs for bifurcating solutions in structural mechanics is \cite{PichiRozza}.

In standard implementations of ROMs, including the papers cited above related to bifurcations, a single {\em global} basis is used to determine the ROM approximation $\uuul$ at any chosen parameter point $\bthe\in\bThe$ by solving \eqref{abs-romg}. However, in a setting in which $\bThe$ consists of subregions for which the corresponding solutions of \eqref{abs-pro} have different character, such as is the case in the bifurcation setting, it may be the case that $\LLL$, although small compared to $\NNN$, may be large enough so that solving the dense ROM system \eqref{abs-romg} many times becomes a costly proposition. Thus, it seems prudent to construct several {\em local} bases, each of which is used for parameters belonging to a different subregion of the bifurcation diagram and also possibly to bridge across the boundary between those subregions. Thus, our goal is to construct, say, $\KKK$ such local bases\footnote{Note that by ``local basis'' we do not mean local with respect to the support of the basis functions in the spatial domain. In fact, the ROM basis function we construct are global in that respect. By ``local basis'' we mean a basis that is used to obtain, using \eqref{abs-roml}, ROM solutions for parameters that lie in subregions of the parameter domain.} of dimension $\LLLk$, each spanning a local subspace $\VVVlk\subset\VVVh$. We then construct $\KKK$ {\em local reduced-order models} 
\be\label{abs-roml}
\ba
    \cN(\uuulk,\vvv;\bthe) = \cF(\vvv;\bthe)\quad\forall\, \vvv\in\VVVlk &\qquad \mbox{for $\kkk=1,\ldots,\KKK$}
    \\ &\Longleftarrow \mbox{\em local reduced-order models}
\ea
\ee
that provide  acceptably accurate approximations $\uuulk$ to the solution $u$ of \eqref{abs-pro} for parameters $\bthe$ belonging to different disjoint parts of the bifurcation diagram. Of course, if one is to employ several local bases, then one must also determine when one needs to switch from one basis to another. 

The use of local ROM bases has been considered in previous works. In \cite{Rapun} snapshots are obtained by sampling the full-order solutions at various time instants; snapshots corresponding to different parameter values are not considered. In \cite{Amsallem1,Amsallem2,Amsallem3} local bases are determined by projection based clustering instead of Euclidean closeness. In 
\cite{kaiser1,kaiser2}, k-means clustering and nearest neighbor classifier with respect to parameters or a low-dimensional representation of the current state are used. In \cite{Peherstorfer} k-means clustering is used to generate local bases for use in conjunction with the empirical interpolation method. In \cite{Zhan} local bases determined using k-means clustering and logistic regression classifiers are used for aero-icing problems. In \cite{Pagani} time-based snapshot clustering as well as k-means and a bisection process clusterings of parameter based snapshots are considered and applied to a model problem in cardiac electrophysiology. In \cite{Ohlberger1,Ohlberger2,Ohlberger3} localized bases are used to resolve fine-scale phenomena in a multiscale setting. In \cite{Eftang} an hp-ROM approach is used to localize ROMs for improved accuracy. In \cite{ChaconRebollo2018} local projection spaces are used in multi-scale turbulence models. None of these efforts specifically address the goal of this paper, i.e., the use of localized ROM basis for bifurcation problems.

The plan for the rest of the paper is as follows. In \S\ref{sec:new}, we discuss in detail the various ingredients that, together, constitute our recipe for the construction and application of a ROM that is well suited for bifurcation problems. Specifically, we discuss how to 

\begin{enumerate}
\item[\bf 1.] {\em select sample points $\{\bthes\}_{\sss=1}^\SSS$ in the parameter domain $\bThe$};

\item[\bf 2.] {\em compute the corresponding {\em snapshots} $\{\uuu_{\NNN,\sss}(\bx)\}_{\sss=1}^\SSS$}, i.e., solutions of the full-order discretized problem \eqref{abs-dis} for each of the parameter points $\bthes$;

\item[\bf 3.] {\em cluster the snapshots so that each cluster corresponds to parameters in a different part of the parameter domain} that could be a subregion of the bifurcation diagram or that could bridge across the boundary of two such subregions;

\item[\bf 4.] {\em construct the local bases corresponding to each cluster};

\item[\bf 5.] {\em detect which cluster a new parameter choice belongs to} so that one can use the correponding local basis to determine a ROM approximate solution.

\end{enumerate}

\noindent When these steps are completed, we can solve \eqref{abs-roml} for the ROM approximation using the appropriate local basis.

In \S\ref{sec:nse}, to illustrate the implementation and employment of the new algorithm defined in \S\ref{sec:new} in a concrete setting, we consider two problems for the Navier-Stokes equations. Concluding remarks are provided in \S\ref{sec:con}.

\section{Detailed description of the new localized-basis method}\label{sec:new}

In this section we provide a more detailed description of each of the five steps listed in the recipe given in \S\ref{sec:intro}.

\subsection{Parameter sampling for snapshot generation}\label{sec:sam}

The recipe outlined in \S\ref{sec:intro} requires, at the start, the selection of $\SSS$ parameter points $\{\bthe_\sss\}_{\sss=1}^\SSS$ in the parameter domain $\bThe\in\R^\MMM$, where $\MMM$ denotes the number of parameters. There exists a large body of literature on methods for this purpose, especially for sampling in hypercubes. Popular sampling techniques include Monte Carlo, quasi-Monte Carlo, Latin hypercube, lattice, orthogonal array, sparse grid, etc. methods. As already mentioned, in this paper we focus on the construction of local reduced bases so that we are content to use crude hand-sampling strategies that pack points near the boundaries between the subregions in the bifurcation diagram. In \S\ref{sec:con}, we briefly comment on more sophisticated sampling strategies that potentially result in more accurate ROMs and which we consider in greater detail in our follow-up paper. 

\subsection{Generating snapshots}\label{sec:snap}

At this point we have in hand a set of $\SSS$ points $\bThes=\{\bthe_\sss\}_{\sss=1}^\SSS$ in the parameter domain $\bThe\subset \R^\MMM$. The second task in the recipe given in \S\ref{sec:intro} is to determine a solution of the full-order model \eqref{abs-dis} for each $\bthe_\sss$, $\sss=1,\ldots,\SSS$. The solution corresponding to $\bthe_\sss$ is denoted by $\uuuhs$ and is referred as the {\em snapshot} corresponding to that parameter point and the set $\cS =\{\uuuhs \}_{\sss=1}^\SSS$ is referred to as the {\em snapshot set}.\footnote{If the problem \eqref{abs-pro} is time dependent, then, for each parameter vector sample $\bthe_\sss$, one can also include in the snapshot set the solution of \eqref{abs-dis} evaluated at several chosen instants of time. Here, to present our ideas, it suffices to treat steady-state problems so, of course, we do not include snapshots in time.}

If $\{\phph_\nnn(\bx)\}_{\nnn=1}^\NNN$ denotes a basis for $\VVVh$, we have, for each $\sss=1,\ldots,\SSS$, that 
\be\label{snaps}
  \uuuhs(\bx)=\sum_{\nnn=1}^\NNN \UUU_{\nnn,\sss} \phph_\nnn(\bx),
\ee
where the coefficients $\{\UUU_{\nnn,\sss}\}_{\nnn=1}^\NNN$ are determined by solving the system of $\NNN$ discrete equations \eqref{abs-dis} for each of the $\SSS$ parameter points $\bthes$, i.e., for each $\sss=1,\ldots,\SSS$, we solve
\be\label{abs-dis2}
    \cN\Big(\sum_{\nnn=1}^\NNN U_{\nnn,\sss} \phph_\nnn(\bx),\phph_{\nnn'};\bthe_\sss\Big) = \cF(\phph_{\nnn'};\bthe_\sss)\qquad\mbox{for $\nnn'=1,\ldots,\NNN$}.
\ee
In cases of greatest interest, the PDE problem \eqref{abs-pro} is nonlinear so that the discrete system \eqref{abs-dis2} is nonlinear as well.

What remains to be specified is the choice of approximating space $\VVVh$, the choice of basis $\{\phph_\nnn(\bx)\}_{\nnn=1}^\NNN$ for $\VVV_\NNN$, and the choice of solution technique used for the nonlinear discrete system of equations \eqref{abs-dis2}. These three choices are usually problem dependent. For the numerical results reported in \S\ref{sec:cha} and \S\ref{sec:cav}, we consider a finite element method \cite{G89,B-quarteroniv2,B-brenners} and a spectral element method \cite{karniadakis1999spectral,CHQZ1,CHQZ2} to effect spatial discretization.

Once the approximating space $\VVVh$ and a basis $\{\phph_\nnn(\bx)\}_{\nnn=1}^\NNN$ for that space are chosen, a snapshot \eqref{snaps} is completely determined by the vector of its coefficients $\u_{\NNN,\sss} = (\UUU_{1,\sss},\ldots,\UUU_{\NNN,\sss})^\top \in \R^\NNN$; naturally, we refer to $\u_{\NNN,\sss}$ as the {\em snapshot vector} corresponding to the snapshot function $\uuuhs(\bx)$ and the set $\cU_{\NNN,\SSS}=\{ \u_{\NNN,\sss} \}_{\sss=1}^\SSS$ as the {\em set of snapshot vectors.} 

It is important to emphasize that the snapshot set is the only information that is available when constructing the ROM we consider using, as is true for most other ROMs. Thus, if the snapshot set does not contain sufficient information to accurately capture desirable features of the solution of the discrete system \eqref{abs-dis} for all parameter values, then the ROM will not be able to do so either. As is said colloquially, ``if it ain't the snapshot set, it ain't in the reduced-order solution.'' On the other hand, the snapshot set itself is completely dependent on the  choice of sample points $\{\bthe_\sss\}_{\sss=1}^\SSS$.

\subsection{Clustering of snapshots}\label{sec:clu}

The third step in the recipe given in \S\ref{sec:intro} is to cluster the snapshots. It is tempting to believe and to therefore require of any clustering algorithm that any cluster of snapshots should correspond to parameter points belonging to a single subregion in the bifurcation diagram. However, bifurcations come in two guises and only for one guise is such an expectation warranted. In the other, one should not expect or desire that a clustering algorithm perform in this manner.

First, there is the case in which passage of the parameter point from one subregion of the bifurcation diagram to another causes {\em discontinuous changes} in some characteristic of the solution of \eqref{abs-pro} or \eqref{abs-dis}, as would happen if one is studying the buckling of a bar under a compressive load; the example in \S\ref{sec:cav} falls into this category as well. The reasonable expectation in this case is that the snapshots  corresponding to parameters in one subregion of the bifurcation diagram are much less ``like'' the snapshots in other subregions than they are to those in its own subregion. Thus, one would expect the members of any cluster of snapshots to correspond to parameters in a single subregion of the bifurcation diagrams, i.e., that {\em clustering the snapshots according to their similarity and clustering parameter points according to which subregion of the bifurcation diagram they belong to would yield the same results.} Note that it may happen and it may even be desirable to happen that more than one snapshot cluster corresponds to the same subregion of the bifurcation diagram.

On the other hand, the dependence of the solution on the parameters could be {\em continuous}, even when passing from one subregion in the bifurcation diagram to another, as would happen if one is studying the stretching of a bar under a tensile load; the example in \S\ref{sec:cha} falls into this category as well. In this case, the clustering situation is not so clear. One would now in fact reasonably expect that the solution and therefore the snapshots near the boundary of a subregion are more ``like'' those on the opposite side of that boundary (because the parameter values are close to each other and because of the continuous dependence on the parameters) compared to snapshots in the same subregion but corresponding to parameter values not close to each other. Thus, in this case, a snapshot cluster could correspond to parameter points in more than one subregion in the bifurcation diagram, i.e., {\em clustering the snapshots according to their similarity and clustering parameter points according to which subregion of the bifurcation diagram they belong to would yield different results.} In the present context, it makes more sense to cluster according to the similarity between snapshots because the ROM bases are built from snapshots so that they influence their construction more directly than the parameters.

\subsubsection{k-means-based clustering of snapshots}
To effect the clustering of the snapshots we use to construct the ROM, we proceed through the following steps, at the heart of which is the use of k-means clustering.
\begin{itemize}
\item {\em For $\kkk=2,3,\ldots$, cluster the set of $\SSS$ snapshot vectors $\cU_{\NNN,\sss}=\{ \u_{\NNN,\sss}\}_{\sss=1}^\SSS$ into $k$ clusters $\{\cV_{\kkk'}\}_{\kkk'=1}^\kkk$ using the k-means algorithm.}
\item {\em Determine the corresponding cluster means $\{\z_{\kkk'}\}_{\kkk'=1}^\kkk$.} 
\item {\em For each $\kkk=2,3,\ldots$, compute the k-means variance  defined in \S\ref{sec:cvt3}.}
\item {\em Stop when a plateau is detected in the k-means variance.}
\end{itemize}

\noindent Note that $\kkk=1$ is omitted from the iterations because, of course, it corresponds to the entire set of snapshots so that no clustering is needed there is simply only the single mean of all the snapshots.

\subsubsection{Determining the optimal number of clusters}\label{sec:cvt3}

A natural question to ask about the clustering of snapshots is how to choose $\KKK$, the number of clusters? For the same input data, i.e., for the same set of snapshots, the k-means variance 
\be\label{cvt-ene2}
{\bm V}(\z_1,\ldots,\z_\KKK; \mathcal{V}_1,\ldots,\mathcal{V}_\KKK) = \sum_{\kkk=1}^\KKK \, \sum_{\{\sss\,:\,\u_{\NNN,\sss}\in\mathcal{V}_\kkk \}}  | \u_{\NNN,\sss}  - \z_\kkk |^2.
\ee
decreases as the number of clusters increases. This is obvious because if we have a k-means clustering with $\KKK$ clusters, then adding another cluster mean at one of the data points $\u_{\NNN,\sss}$ would clearly reduce the variance. Furthermore, if $\KKK=\SSS$, then we may pick all the generators to coincide with the data points, leaving us with zero variance but with each cluster containing a single point. Clearly, that is not what one wants, so the question posed above needs a more rational answers.

\vskip5pt
\noindent{\textbf{\em Variance elbow.}}
If the input data is amenable to clustering, then, for a small number of clusters, the k-means variance reduces quickly as the number of clusters increases. However, eventually the decrease in the variance lessens as one adds more clusters. As an example, consider the discrete data case for which the $\SSS$ input vectors are distributed evenly among four balls of radius 1 but whose centers form a square of side 10. Clearly, going from 2 to 3 k-means clusters would reduce the k-means variance substantially. This is again true when going from 3 to 4 clusters. However, clearly, going from 4 to 5 clusters does not buy one as much reduction in the variance because of the nature of the data.

Thus, if the input data is reasonably clustered, one can expect the plot of the k-means variance vs. the number of clusters to have an ``elbow,'' i.e., at a certain number of clusters, the plot transitions from having a steep negative slope for a low number of clusters to a gently sloping plateau for a higher number of clusters. Quantitatively, one can use the elbowing effect to choose the number of clusters $\KKK$ to be the smallest integer for which the difference in the k-means variances for $\KKK-1$ and $\KKK$ clusters is a specified fraction of the difference in the k-means variance between 2 and 3 clusters.

In \S\ref{sec:cha} and \S\ref{sec:cav}, we provide plots of the k-means variance for the examples we consider, plots which exhibit the elbow shape.

\subsection{Construction of local ROM bases}\label{sec:bas}

The fourth step in the recipe given in \S\ref{sec:intro} is to construct local ROM bases, i.e., low-dimensional bases for use in \eqref{abs-roml} to obtain ROM approximations of the solution of \eqref{abs-pro}, i.e., each local ROM basis is meant to help account for a separate subregion of the parameter domain $\bThe$, including those that straddle across the boundaries between two subregions of the bifurcation diagram. 

Perhaps the two approaches in most common use for the construction of ROM bases are the {\em reduced-basis}\footnote{All ROMs of the type we consider in this paper involve ``reduced'' bases in the sense that they involve very low-dimensional bases compared to the bases used to effect spatial approximations. However, the nomenclature ``reduced basis method'' is generally reserved for the context of using all the snapshots as the ROM basis.} (see, e.g., \cite{hesthaven2015certified, QuarteroniManzoniNegri}) and {\em proper orthogonal decomposition} (see, e.g, \cite{Volkwein}) methods. RB methods usually use the whole snapshot set as the ROM basis whereas POD bases are determined through an optimization process that selects a basis that spans a subset of the span of the snapshot vectors. A third alternative is to use the cluster means themselves as a ROM basis. Here, so as to have a concrete setting for our numerical illustrations, we choose to use POD bases.

We have in hand the set of snapshot vectors $\cU_\SSS = \{ \u_{\NNN,\sss} \}_{\sss=1}^\SSS$ and a clustering $\cU_\SSS=\cup_{\kkk=1}^\KKK \cV_\kkk$
of that set into disjoint subsets $\cV_\kkk$, $\kkk=1,\ldots,\KKK$, with $\cvk$ denoting the cardinality of $\cV_\kkk$; see \S\ref{sec:clu}. Using this data, we define $\KKK$ {\em local POD bases} $\{ \ppsi_{\lll,\kkk} \}_{\lll=1}^{\LLL_\kkk}$ with $\LLL_\kkk\le\cvk$, one for each cluster $\cV_\kkk$, by solving, for $k=1,\ldots,K$, the Euclidean norm minimization problem
\be\label{pod:loc}
\ba
&\min_{\ppsi_{1,\kkk},\ldots,\ppsi_{\LLL_\kkk,\kkk} }
\sum_{\{\sss=1\,\,:\,\,\u_{\NNN,\sss}\in\cV_\kkk\}} 
\Big|
\u_{\NNN,\sss}-\sum_{\ell=1}^{\LLL_\kkk} 
\big(\u_{\NNN,\sss} \cdot\ppsi_{\lll,\kkk}\big)  \ppsi_{\lll,\kkk}
\Big|^2 
\\&\qquad\qquad\qquad
\mbox{such that $\ppsi_{\lll,\kkk}\cdot\ppsi_{\lll,\kkk'}=\delta_{\kkk\kkk'}$ for $\kkk,\kkk'=1,\ldots,\LLL_\KKK$},
\ea
\ee
where $\delta_{\kkk\kkk'}$ denotes the Kroenecker delta. If $\KKK=1$ there is, of course, no clustering so that the single POD basis is based on all $\SSS$ snapshots; in this case, we refer to the basis as the {\em global POD basis.}

For each cluster $\cV_\kkk$, $\kkk=1,\ldots,\KKK$, of snapshots, we define the {\em cluster snapshot matrix} $\SM_\kkk$ as the $\NNN\times\cvk$ matrix whose columns are the $\NNN$-dimensional snapshots vectors in cluster $\cV_\kkk$. Then, for each $\kkk$, the solution of the minimization problem \eqref{pod:loc} can be obtained from the singular value decomposition (SVD) $\SM_\kkk =  \Psps_\kkk  \bSigma_\kkk  \Upup^\top$ of the cluster snapshot matrix. Here, $\Psps_\kkk$ and $\Upup_\kkk$ are orthogonal $\NNN\times\NNN$ and $\cvk\times\cvk$ matrices, respectively, and $\bSigma_\kkk$ is an $\NNN\times\cvk$ matrix with all its nonzero entries appearing along the main diagonal in non-increasing order, i.e., $\sigma_{1,\kkk}\ge\sigma_{2,\kkk}\ge\cdots\ge\sigma_{\cvk,\kkk}\ge0$, where $\{\sigma_{\lll,\kkk}\}_{\lll=1}^\cvk$ denote the singular values of $\SM_\kkk$. The columns of $\Psps_\kkk$ and $\Upup_\kkk$ are referred to as the {\em left and right singular vectors} of $\SM_\kkk$, respectively, and the nonzero entries of $\bSigma_\kkk$ are referred to as the {\em singular values} of $\SM_\kkk$. For each $\kkk$, the solution of the minimization problem \eqref{pod:loc} is given by the first $\LLL_\kkk$ columns of $\Psps_\kkk$, i.e., {\em the $\LLL_\kkk$ local POD basis vectors corresponding to the cluster $\cV_\kkk$ is given by $\{\ppsi_{\lll,\kkk} \}_{\lll=1}^{\LLL_k}$, where $\ppsi_{\lll,\kkk}$ is the $\lll$-th left-singular vector of $\SM_\kkk$.}

Of course, each POD basis vector induce a function in the space $\VVVh$ used for the spatial discretization \eqref{abs-dis} of the given problem \eqref{abs-pro}. Again letting $\{\phph_\nnn(\bx)\}_{\nnn=1}^\NNN$ denote a basis for $\VVVh$, we have, corresponding to the cluster $\cV_\kkk$, the POD basis functions
\be\label{basisl}
  \vvv_{\lll,\kkk}(\bx)   = \sum_{\nnn=1}^\NNN \psi_{\nnn,\lll,\kkk}  \phph_\nnn(\bx) \qquad\mbox{for $\lll=1,\ldots,\LLLk$},
\ee
where $\psps_{\nnn,\lll,\kkk}$ denotes the $n$-th component of the POD basis vector $\ppsi_{\lll,\kkk}$. Then, the subspace of the spatial approximation space $\VVVh$ corresponding to the cluster $\cV_\kkk$ in snapshot space is given by $\VVV_{\LLLk}={\rm span}\{\vvv_{\lll,\kkk}\}_{\lll=1}^\LLLk$. Note that the error in \eqref{pod:loc} is given by the sum of the neglected singular values $\sigma_{\lll,\kkk}$ for $\lll=\LLL_\kkk+1,\ldots, \cvk$. Hence, one might choose the number of basis functions $\LLL_\kkk$ such that $\sum_{\lll=\LLL_\kkk+1}^{\cvk} \sigma_{\lll,\kkk}$ is smaller than some prescribed tolerance. We further note that this is a heuristic consideration because the problem \eqref{pod:loc} is dependent on the snapshots.

\subsection{Determining which local basis to use for a new parameter}\label{sec:bif_detection}

In the fifth step in the recipe given in \S\ref{sec:intro}, the task faced is to identify, for any parameter point $\bthe\in\bThe$ which was not among those used to generate the snapshots, which local basis should be used in \eqref{abs-roml} to determine the corresponding ROM approximation. The local bases are associated with the clustering of snapshots so seemingly, to identify which local basis one should use, one should solve the expensive full-order model \eqref{abs-dis} for that $\bthe$ and then determine which of the $\KKK$ cluster means of the snapshot clusters $\{\cV_\kkk\}_{\kkk=1}^\KKK$ that solution is closest to. Obviously one does not want to do this because involves costs depending on the number of full-order degrees of freedom. Thus, {\em although the snapshots are used to determine the clusters, the identification of the local basis to be used has to be done in parameter space.} 

Thus the first step in any recipe for assigning new parameters to clusters is to {\em determine the parameter clustering induced by the snapshot clustering.} There is, of course, a one-to-one correspondence between a snapshot vector $\u_{\NNN,\sss}$ (or the corresponding snapshot function $\uuuhs(\bx)$) and the parameter point $\bthes$ used to determine that snapshot. Thus, a clustering $\{\cV_\kkk\}_{\kkk=1}^\KKK$ of the snapshots trivially induces a clustering of the sample points $\bThes=\{\bthes\}_{\sss=1}^\SSS$. Specifically, we have that $\bThes=\cup_{\kkk=1}^\KKK \bThe_\kkk$, where, for $\kkk=1,\ldots,\KKK$, 
\be\label{parclu}
   \bThe_\kkk = \{ \bthes\in\bThes \,\,:\,\,
   \u_{\NNN,\sss}\in\cV_\kkk\}.
\ee

There are potential difficulties stemming from the fact that the clustering is done with respect to the snapshots whereas the identification of which basis to use is done with respect to the parameters. For example, a tempting way to assign a new parameter to a cluster (and thus select which local basis to use for solving the ROM model for that parameter) is to first determine the cluster means in parameter space and then assign the new parameter to the cluster corresponding to the closest cluster mean. However, a parameter cluster mean generally does not correspond to the snapshot cluster mean, i.e., if one were to solve the full-order model using the parameter cluster mean, the resulting full-order solution would be different from the snapshot cluster mean. As a result, the possibility exists that the wrong local ROM basis would be assigned to a new parameter. This is much more likely to occur in the case where transitions across boundaries in the bifurcation diagram result in continuous solutions and it is likely that clusters span across such boundaries. In this case, because
the transition results in continuous solutions, the ROM approximation error when using the ``wrong'' local basis still be tolerable, i.e., less than $1\%$. In the rarer case in which the transition across boundaries in the bifurcation diagram is discontinuous, using the ``wrong'' local basis can e expected to result in relatively large approximation error.

In one dimension, i.e., for the case of a single parameter, an alternate means for assigning new parameters to parameter clusters makes use of the cluster midrange and cluster radius and proceeds as follows. We have in hand the clustering $\{\cV_\kkk\}_{\kkk=1}^\KKK$ of the snapshots determined as described in \S\ref{sec:clu} and a point $\bthe\in\bThe$ for which we want to obtain the corresponding ROM solution. Then, we
\begin{itemize}
\item determine, using \eqref{parclu}, the clustering $\{\bThe_\kkk\}_{\kkk=1}^\KKK$ in the parameter domain induced by the clustering of snapshots;

\item compute, for $\kkk=1,\ldots,\KKK$, the midrange 
$$
   {M}^{(mr)}_\kkk= \frac1{2} \Big( \min_{\{\sss=1,\ldots,\SSS,\,\,\,\theta_s\in\,\bThe_\kkk\}} \theta_s + \max_{\{\sss=1,\ldots,\SSS,\,\,\,\theta_s\in\,\bThe_\kkk\}} \theta_s \Big) 
$$
and the radius 
$$
   {\overline{r}}_\kkk= {M}^{(mr)}_\kkk - \min_{\{\sss=1,\ldots,\SSS,\,\,\,\theta_s\in\,\bThe_\kkk\}} \theta_s  
$$
of the parameters in each parameter cluster;

\item assign the new parameter to the parameter cluster $\bThe_\kkk$ such that
$$
   |\theta - M^{(mr)}_\kkk| - {\overline{r}}_\kkk < |\theta - M^{(mr)}_{\kkk'}| - {\overline{r}}_{\kkk'} \qquad
   \mbox{for $\kkk'=1,\ldots,\KKK,\,\,\,\kkk'\ne\kkk$},
$$
i.e, the parameter cluster whose midrange minus the cluster radius is closest to that parameter, note that the term $|\bthe - {\overline{m}}_\kkk| - {\overline{r}}_\kkk$ can also be negative;

\item the parameter cluster $\bThe_\kkk$ so identified corresponds to a specific snapshot cluster $\cV_\kkk$;

\item the local basis to use is $\{\ppsi_{\lll,\kkk} \}_{\lll=1}^{\LLL_k}$, i.e., the local basis corresponding to snapshot cluster $\cV_\kkk$.
\end{itemize}

\noindent Limited computational experiments indicate that this may be a superior avenue for assigning a local basis to a new parameter when compared to using parameter cluster means. Of course, there is no shortage of ways to generalize the concepts of midrange and radius to higher dimensions, i.e., to the case of multiple parameters. For example, for the midrange, one could choose the center of the smallest ball containing all the parameters within a cluster, or one could average the positions of the vertices of the convex hull polytope of the parameters, or one could determine the midrange for each component of the parameter vector and use the vector of midranges as the testing point. The corresponding radii would be the radius of the ball, the distance from the midrange to the furthest vertex, and the largest of the one-dimensional radii. In one dimension, all three produce the midrange point.

\section{Application to the incompressible Navier-Stokes equations}
\label{sec:nse}

We use a variational formulation of the incompressible Navier-Stokes equations as the concrete setting of \eqref{abs-pro} to illustrate our methodology. We first define the Navier-Stokes problem and then consider examples for both continuous and discontinuous transitions through bifurcation points. For the first case, we consider two types of spatial discretizations. Together, the examples are meant to illustrate the robustness of our approach with respect to both solution behaviors and spatial discretization choices.


\subsection{Problem definition}\label{sec:pd}

The Navier-Stokes equations for incompressible flows are given by
\begin{equation}\label{NS-1}
\begin{aligned}
\frac{\partial {\bm u}}{\partial t}+({\bm u}\cdot  \nabla {\bm u})-\nu\Delta {\bm u}+\nabla p={\bm f} &\qquad\mbox{in } \Omega\times(0,T]\\
\nabla \cdot {\bm u}=0&\qquad\mbox{in } \Omega\times (0,T],
\end{aligned}
\end{equation}
where ${\bf u}$ and $p$ denote the unknown velocity and pressure fields, respectively, $\nu>0$ denotes the kinematic viscosity of the fluid, and ${\bm f}$ accounts for possible body forces. These equations describe the incompressible motion of a viscous, Newtonian fluid in the spatial domain $\Omega \subset\mathbb{R}^d$, $d = 2$ or $3$, over a time interval of interest $(0, T]$. The boundary of $\Omega$ is denoted by $\partial\Omega$. If the fluid acceleration becomes negligible, i.e., the system has evolved towards a steady state, the time-derivative term in \eqref{NS-1} can be omitted.

To complete the definition of the flow problem, \eqref{NS-1} needs to be endowed with initial and boundary conditions, which, for the examples of \S\ref{sec:cha} and \S\ref{sec:cav}, take the form
\begin{eqnarray}
{\bm u}={\bm u}_0&\qquad&\mbox{in } \Omega \times \{0\} \label{IC}\\
{\bm u}={\bm u}_D&\qquad&\mbox{on } \partial\Omega_D\times (0,T] \label{BC-D} \\
-p \n + \nu \frac{\partial {\bm u}}{\partial \n}=\g &\qquad&\mbox{on } \partial\Omega_N\times, (0,T], \label{BC-N}
\end{eqnarray}
where $\partial\Omega_D \cap \partial\Omega_N = \emptyset$ and $\overline{\partial\Omega_D} \cup \overline{\partial\Omega_N} =  \overline{\partial\Omega}$. Here, ${\bm u}_0$, ${\bm u}_D$, and $\g$ are given and $\n$ denotes the unit normal vector on the boundary $\partial\Omega_N$ directed outwards. In the rest of this section, we will explicitly denote the dependence of the solution of the problem \eqref{NS-1}-\eqref{BC-N}, and possibly of ${\bm f}$, on the parameter vector $\bthe$.

To define a variational formulation of system~\eqref{NS-1}-\eqref{BC-N}, we let $L^2(\Omega)$ denote the space of square integrable functions in $\Omega$ and $H^1(\Omega)$ the space of functions belonging to $L^2(\Omega)$ with first derivatives in $L^2(\Omega)$. We then define the sets 
$$
{\bm V} := \left\{ {\bm v} \in [H^1(\Omega)]^d: ~ {\bm v} = {\bm u}_D \mbox{ on }\partial\Omega_D \right\}, \quad
{\bm V}_0:=\left\{{\bm v} \in [H^1(\Omega)]^d: ~ {\bm v} = \boldsymbol{0} \mbox{ on }\partial\Omega_D \right\}. \el
$$
Then, the standard variational form corresponding to \eqref{NS-1}-\eqref{BC-N} is given as follows: find $({\bm u}(\bthe),p(\bthe))\in {\bm V} \times L^2(\Omega)$ satisfying the initial condition \eqref{IC} and
\begin{equation}\label{eq:weakNS-1}
\begin{aligned}
&\int_{\Omega} \frac{\partial{\bm u}(\bthe)}{\partial t}\cdot{\bm v}\de\mathbf{x}+\int_{\Omega}\left({\bm u}(\bthe)\cdot\nabla {\bm u}\right)\cdot{\bm v}\de \mathbf{x} - \int_{\Omega}p(\bthe)\nabla \cdot{\bm v}\de\mathbf{x}\\
&\hspace{3cm} =\int_{\Omega}{\bf f}(\bthe)\cdot{\bm v}\de\mathbf{x}  + \int_{\partial \Omega_N}{\bf g}\cdot{\bm v}\de\mathbf{x} 
\qquad\forall\,{\bm v} \in {\bm V}_0 \\
& \int_{\Omega}q\nabla \cdot{\bm u}(\bthe)\de\mathbf{x} =0 \qquad\forall\, q \in L^2(\Omega).  
\end{aligned}
\end{equation}
Problem \eqref{eq:weakNS-1} constitutes the particular case of the abstract problem \eqref{abs-pro} we use for the numerical illustrations.

The nonlinearity in problem \eqref{eq:weakNS-1} can produce a loss of uniqueness of the solution, with multiple solutions branching from a known solution at a bifurcation point. We consider, in \S\ref{sec:cha} and \S\ref{sec:cav}, two specific problems for which the parameter domain $\bThe$ contains two or more regions across whose common boundaries steady bifurcations occur. Because we are interested in studying flow problems close to steady bifurcation points, our snapshot sets include only steady-state solutions as explained in  \S\ref{sec:snap}. To obtain the snapshots, we approximate the solution of problem \eqref{eq:weakNS-1} by a time-marching scheme that we stop when sufficiently close to the steady state, e.g., when the stopping condition
$
\frac{\|{\bf u}_N^n-{\bf u}_N^{n-1}\|_{L^2(\Omega)}}{\|{\bf u}_N^n\|_{L^2(\Omega)}}<\mathtt{tol}
\label{eq:def_increment}
$
is satisfied for a prescribed tolerance $\mathtt{tol}>0$, where $n$ denotes the time-step index.

\subsection{Example 1: Jet emanating from an orifice into a channel} \label{sec:cha}

We first consider a classical benchmark test: a jet exiting through an orifice into a rectangular channel. The simplicity of the geometry that is partially characterized by the expansion ratio, i.e., the ratio of the orifice size to the channel height, but which still results in complex flow have made this problem a popular choice for use in testing computational models; see, e.g., \cite{fearnm1,drikakis1,hawar1,mishraj1} and the references cited therein.

In \eqref{NS-1}, we set $ \f = {\bm 0} $. The channel height is set to $6$. At the orifice, which is of size $w=1$ and is centered on the left side of the channel, we impose \eqref{BC-D} with ${\bm u}_D$ having a parabolic horizontal velocity component with maximum horizontal speed $U=1/4$ and zero vertical component. At the remainder of the left end of the channel and at the top and bottom of the channel, we apply \eqref{BC-D} with ${\bm u}_D = {\bf 0}$. The right end of the channel is an outlet at which \eqref{BC-N} is imposed with $\g={\bf 0}$. With this condition along with the other data imposed, if the channel is sufficiently long, the outflow will be fully established \cite{fortinj1}. For this reason, we choose the channel length to be 36. The Reynolds number $\text{Re} = {U w}/{\nu}$ can be used to characterize the flow regime; it can be thought of as the ratio of inertial forces to viscous forces. For large values of the Reynolds numbers, inertial forces are dominant over viscous forces and vice versa. For the data chosen above, we have that Re $= 1/ (4\nu)$.

For a fixed expansion ratio, as the Reynolds number $\text{Re}$ increases from zero, the sequence of flow configurations is as follows. For sufficiently small values of $\text{Re}$, a steady {\em symmetric} jet is observed. Moffatt eddies form close to the corners downstream of the expansion. As the inertial effects of fluid become more important, the Moffatt eddies develop into two recirculation regions of equal size. As $\text{Re}$ increases further, flow symmetry about the horizontal centerline is initially maintained and the recirculation length progressively increases. However, at the critical value $\text{Re} \approx 33$, one recirculation zone expands whereas the other shrinks, giving rise to a steady {\em asymmetric} jet. This asymmetric solution remains stable as $\text{Re}$ increase further, but the asymmetry becomes more pronounced, as shown in \cite{mishraj1}. The configuration with a symmetric jet is still a solution, but is unstable \cite{sobeyd1}. Snapshots of the flow fields for three stable scenarios are illustrated in  \S\ref{sec322}. Changes in the solution as $\text{Re}$ increases are {\em continuous}, even when passing from a symmetric to an asymmetric configuration. This loss of symmetry in the steady solution as $\text{Re}$ changes is a supercritical pitchfork bifurcation \cite{Prodi}.

\subsubsection{Spatial discretization approaches}\label{sec:jetdisc}

For this example, we compare the use of two spatial discretization methods, specifically, a finite element method (FEM) and a spectral element method (SEM). Our goal here is to show the possible effect that spatial discretization choices can have on the results obtained for our ROM approach. 

\begin{figure}[h!]
\begin{center}
\includegraphics[width=.95\textwidth]{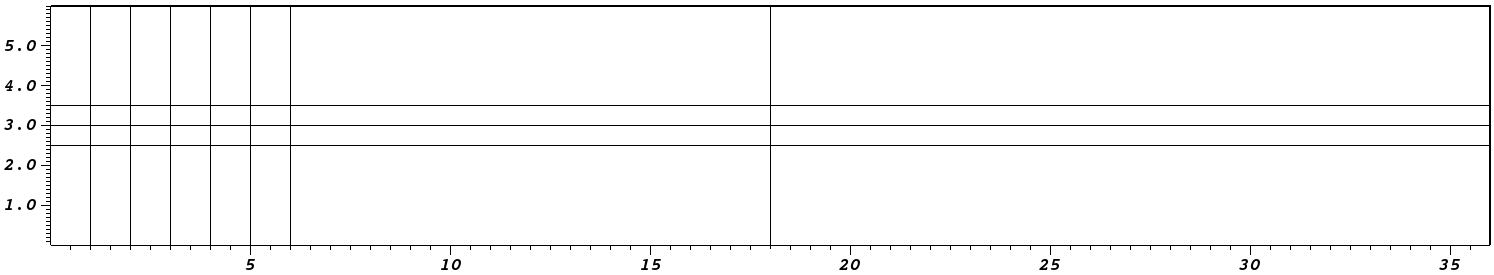}  
\caption{The $32$ spectral elements used for spatial approximation which results in $14259$ degrees of freedom.}
\label{fig_SEM_domain_channel}
\end{center}
\end{figure}

For the {\em spectral element discretization}, the SEM software framework Nektar++, version 4.4.0, (see {\tt https://www.nektar.info/}) is used for the full-order simulations. The domain is discretized into $32$ quadrilateral elements as shown in Fig. \ref{fig_SEM_domain_channel}. Modal Legendre ansatz functions of order $12$ are used in every element and for every solution component. This results in $4753$ global degrees of freedom for each of the horizontal and vertical velocity components and the pressure for the time-dependent simulations. For temporal discretization, an IMEX scheme of order 2 is used with a time-step size of $\Delta t = 10^{-3}$; typically $10^6$ time steps are needed to reach a steady state.

For the {\em finite element discretization,} the FEM software framework FEniCS (see {\tt https://fenicsproject.org/}) is used for the full-order simulations. We effect spatial discretization using the Taylor-Hood finite element pair, i.e., continuous piecewise quadratic polynomials for approximating the velocity components and continuous piecewise linear polynomials for approximating the pressure, all based on the same grid. The triangular grid used results in $90,876$ degrees of freedom. For the temporal discretization, a modified version of Chorin's method, referred to as the incremental pressure correction scheme (IPCS) \cite{G79}, is used with a time step  $\Delta t = 0.1$.  

The numerical study show that a supercritical pitchfork bifurcation occurs around the critical Reynolds number $\approx 33$. 
For values of $\text{Re}$ lower than that critical value, a single branch of symmetric stable approximate solutions are obtained. 
For values above the critical value, three different steady state solutions are possible. 
Besides a solution which is symmetric with respect to the horizontal midline of the domain, two other branches of approximate solutions that are asymmetric with respect to the 
horizontal midline can occur. The symmetric approximate solution is now unstable whereas the asymmetric approximate solutions are stable. This is all in accordance to what is 
known about the exact solution.

\begin{figure}[h!]
  \begin{center}
\includegraphics[width=0.45\textwidth]{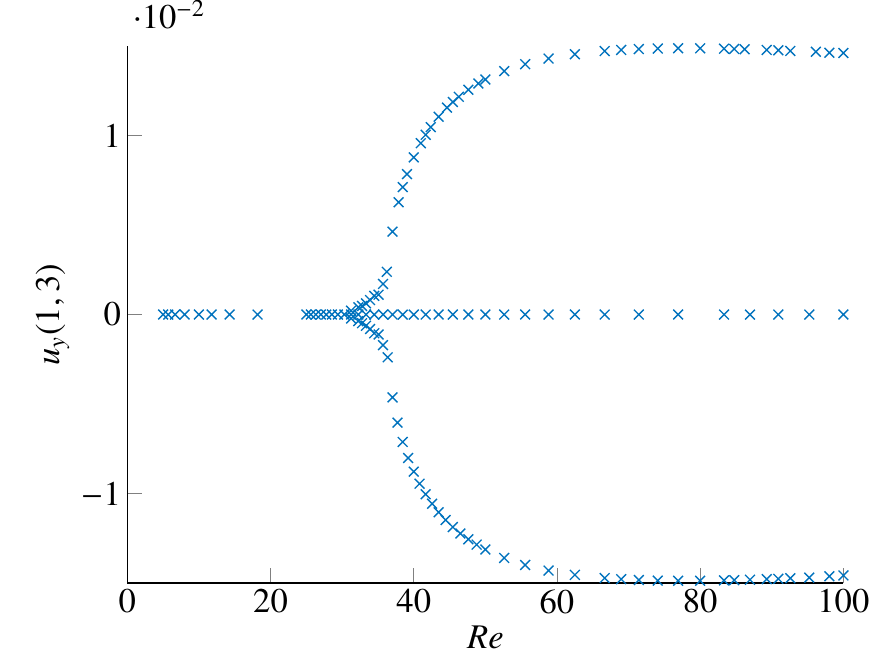}  
    \end{center}   
\caption{For the jet problem, the full bifurcation diagram, based on the vertical velocity component at the point $(1.0, 3.0)$ for Reynolds numbers from $5$ to $100$.} \label{fig:Channel_bifurcation_diagram}   
\end{figure}

To obtain the unstable symmetric approximate solutions for higher Reynolds numbers, a continuation method is used in which the solution is determined with an initial 
condition given from a solution at a slightly lower Reynolds number. To transition onto one of the two stable branches, a solution is first computed with a non-symmetric 
inflow condition, i.e., a stronger inflow towards the either the upper or lower wall, as desired. This solution is then used as an initial condition with a symmetric inflow 
condition, as specified in the model description, to obtain an approximate solution on the desired branch. The desired branch is then followed by again using a continuation 
method to determine approximate solutions for higher values of the Reynolds numbers.  Using theses strategies, the full pitchfork bifurcation diagram can be obtained as illustrated in Fig.~\ref{fig:Channel_bifurcation_diagram} that provides that diagram based on the point value of the approximate vertical velocity component at the point $(1.0, 3.0)$. 

\noindent
When computing the snapshots needed to build a reduced-order model, one clearly needs to follow one or the other of the asymmetric branches i.e., to not mix 
snapshots obtained from both branches. In our calculations, we use the lower branch of asymmetric solutions as determined 
using the continuation strategy mentioned above.

\subsubsection{Results for the offline stage: the construction of local ROM bases}\label{sec322}
\begin{figure}[h!]
\begin{center}
\includegraphics[width=2.25in]{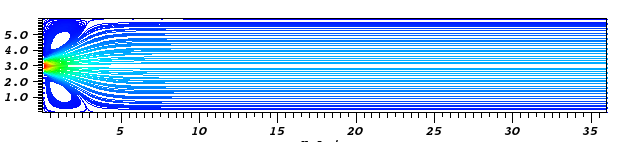} 
\includegraphics[width=2.25in]{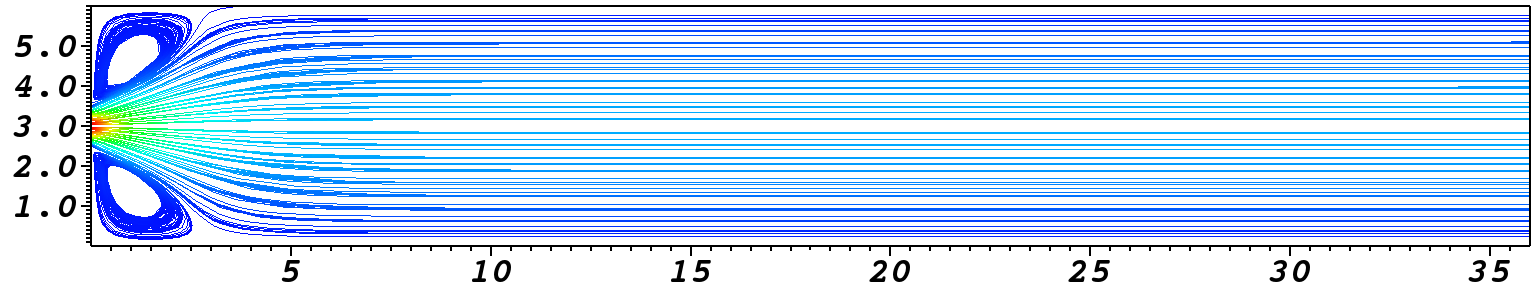}  
\includegraphics[width=.5in]{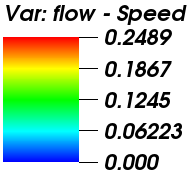}  
 \\
\vskip5pt
 \includegraphics[width=2.25in]{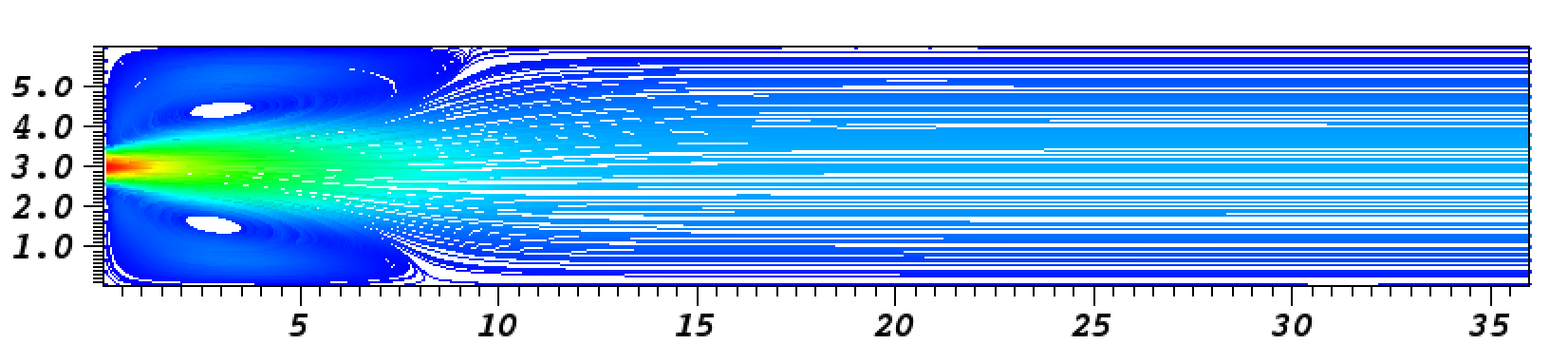}  
\includegraphics[width=2.25in]{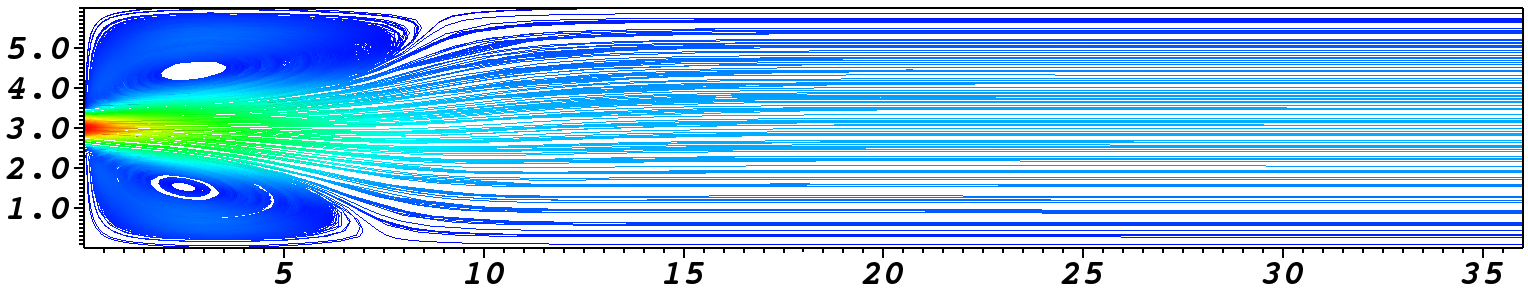}
\includegraphics[width=.5in]{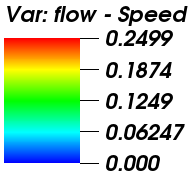} \\
\vskip5pt
\includegraphics[width=2.25in]{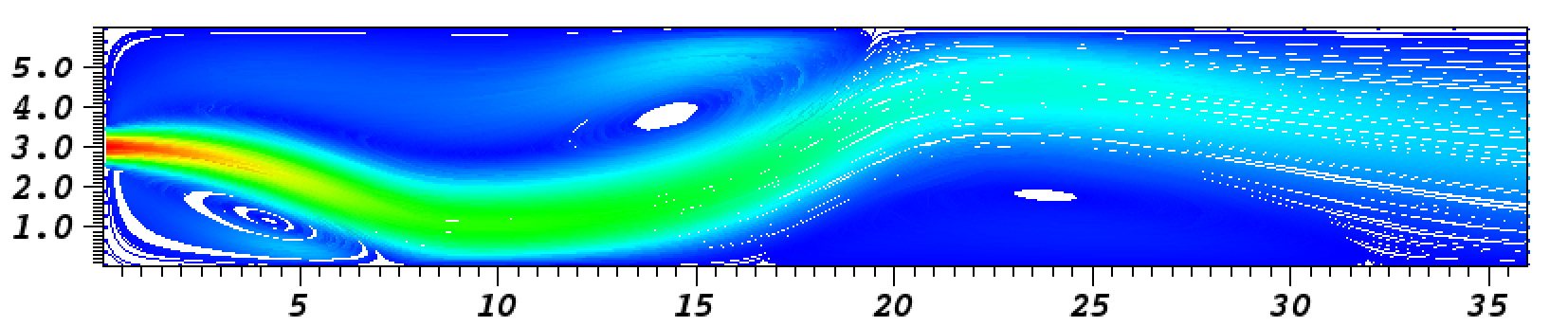} 
\includegraphics[width=2.25in]{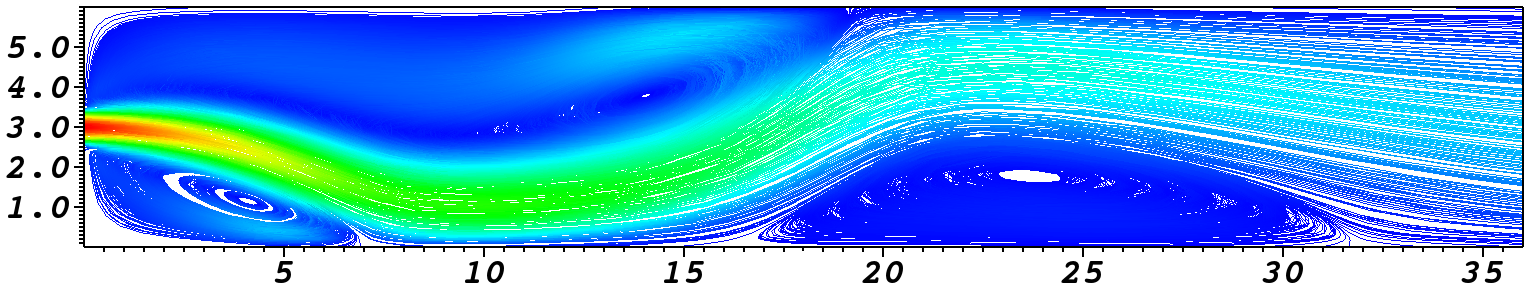}
\includegraphics[width=.5in]{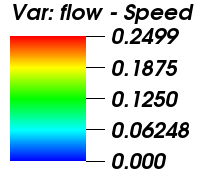}
 \caption{Top to bottom: full-order steady-state streamlines for $\text{Re} =10$, $\text{Re}=35$, and $\text{Re}=100$
 given by FEM (left) and SEM (right). Colors correspond to the magnitude of the velocity according to the given scales.} \label{fig:Channel_full_order}
 \label{fig_FEM_SEM_Re10_35_100stream}
\end{center}
\end{figure}

We next compare results obtained using the two spatial discretization methods for the offline stage in which local ROM bases are constructed.
\vskip5pt
\noindent\textbf{\em Snaphots.}
First, in Fig. \ref{fig_FEM_SEM_Re10_35_100stream}, we provide visual comparisons of steady-state streamline snapshots obtained using the two methods and for three values of the Reynolds number: $\text{Re} =10$ which is below the critical Reynolds number and $\text{Re}=35$ and $100$ which are slightly above and well above, respectively, the critical Reynolds number. Not surprisingly, the results for two methods are in very good agreement, with the visible differences occurring in very low velocity regions which is unavoidable because the discretization errors of the schemes are most probably larger than the speed in those regions.

\begin{figure}[h!]
  \begin{center} 
\includegraphics[width=0.475\textwidth]{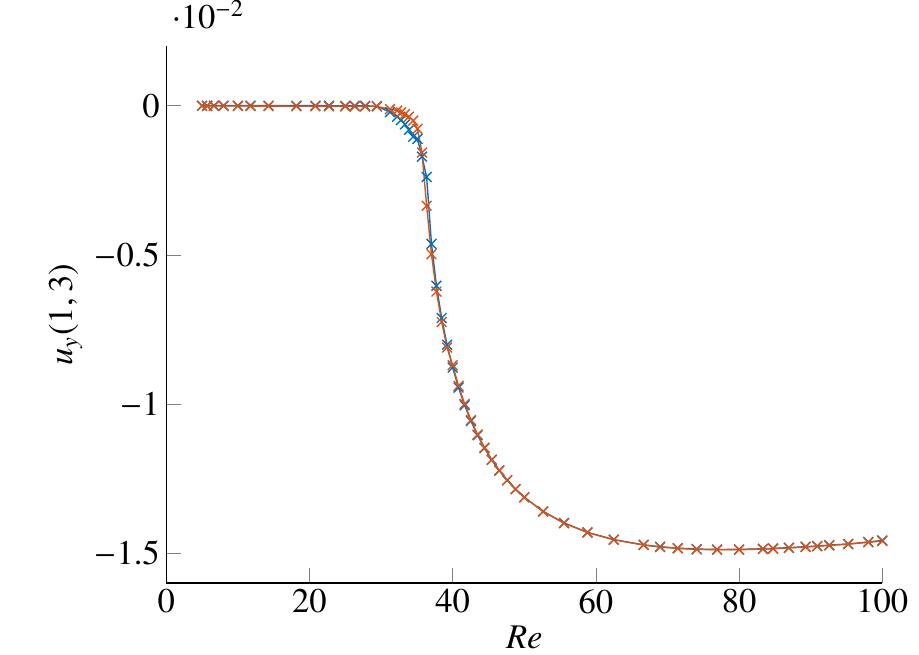} 
    \end{center}   
\caption{Lower stable branch of the bifurcation diagram over Reynolds numbers from $5$ to $100$ for the jet problem for the FEM (red) and SEM (blue) discretizations.} \label{fig:Channel_bifurcation_diagram_lower}
\end{figure}

\vskip5pt
\noindent\textbf{\em Bifurcation diagram.}
We next compare, in Fig. \ref{fig:Channel_bifurcation_diagram_lower}, the (lower half of the) bifurcation diagram for the vertical velocity $u_y$ at the point $(1,3)$ as obtained by the SEM and FEM methods. The two diagrams are in very close agreement, with the only visual differences occurring for Reynolds numbers such that the flow is only very slightly asymmetric and thus the vertical velocity is extremely small. Again, this discrepancy is unavoidable because the discretization errors of the schemes are most probably larger than $|u_y(1,3)|$.

\begin{figure}[h!]
\begin{center}
  \includegraphics[scale=0.7]{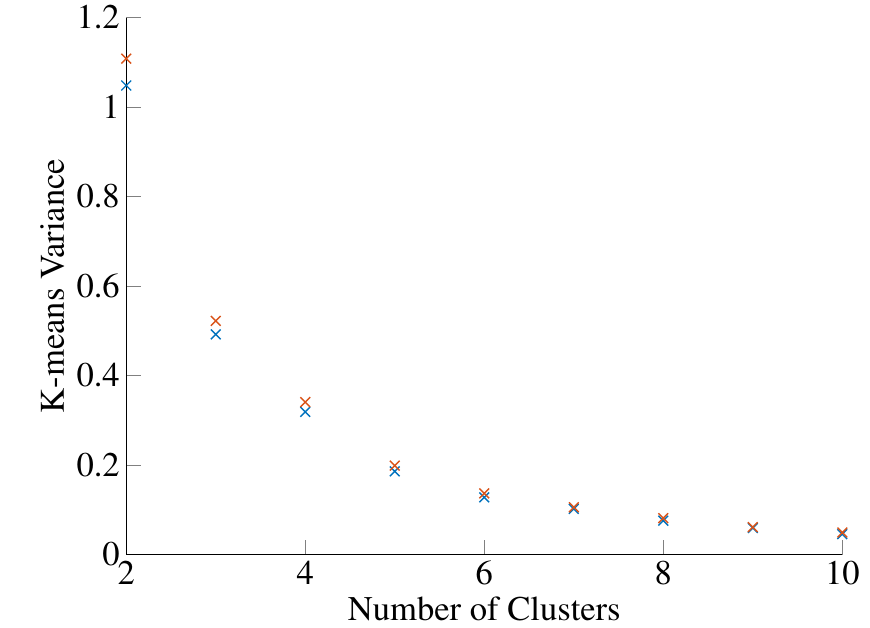}  
\caption{k-means variance versus number of clusters for the jet problem for the SEM (blue crosses) and FEM (red crosses) discretizations.} \label{fig:Channel_CVT}
\end{center}   
\end{figure}

\vskip5pt
\noindent\textbf{\em k-means variance.} The k-means variance is plotted vs. the number of clusters in Fig. \ref{fig:Channel_CVT} for both spatial discretization approaches. We use the elbowing effect of the k-means variance to determine the optimal number of clusters so the close agreement between the SEM and FEM results would result in the same choice for the optimal number of clusters. With regards to the plot itself, it is clear that the value of the k-means variance reduces quickly for a very few clusters and that the reduction lessens as the number of clusters increases. The decision about how many clusters to use should be guided by the goal of making the clustering error commensurate with the other errors incurred in the ROM process.    

\begin{figure}[h!]
\begin{center}
\includegraphics[scale=0.45]{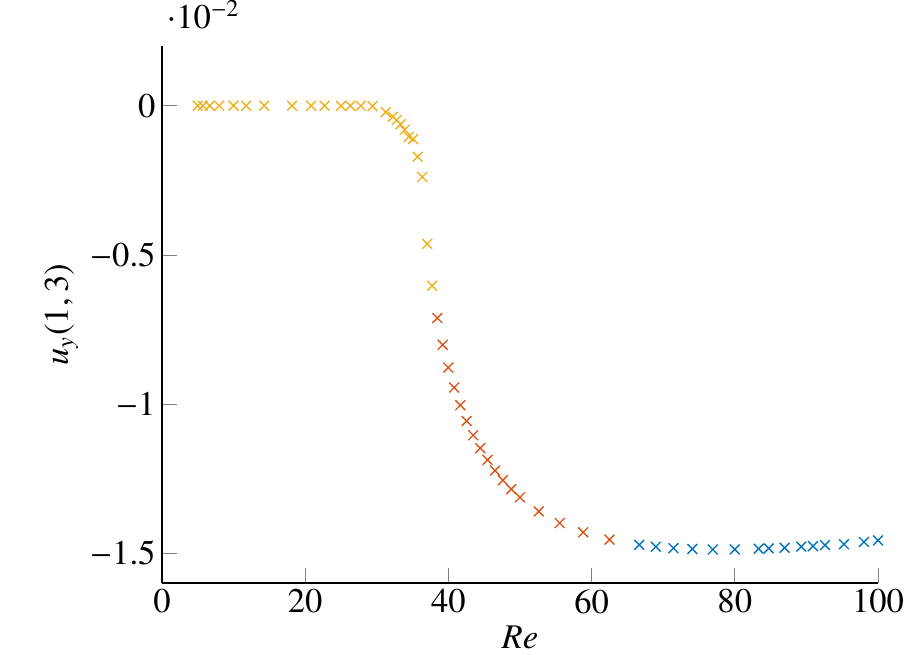}
\includegraphics[scale=0.45]{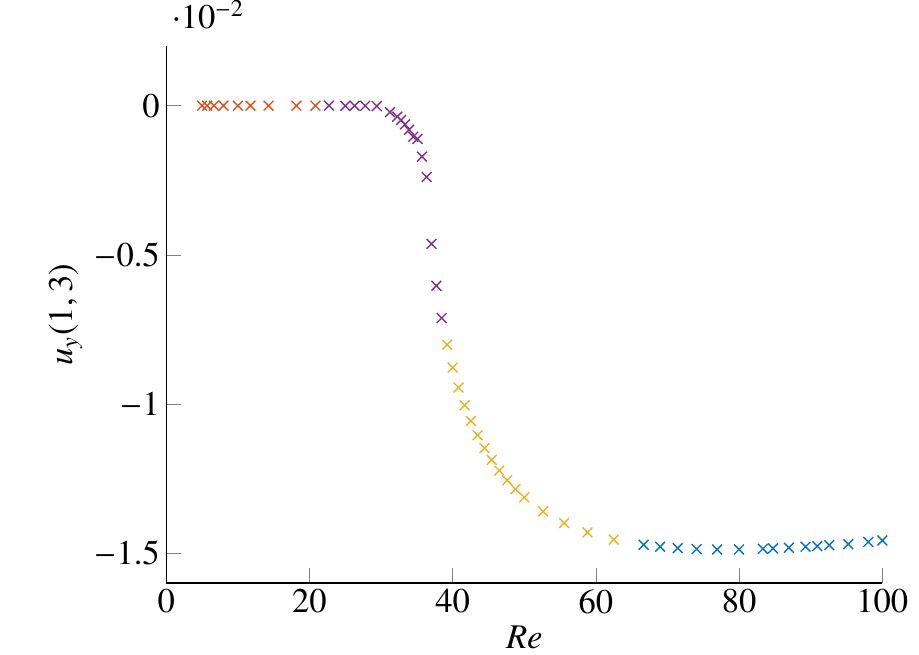}
\includegraphics[scale=0.45]{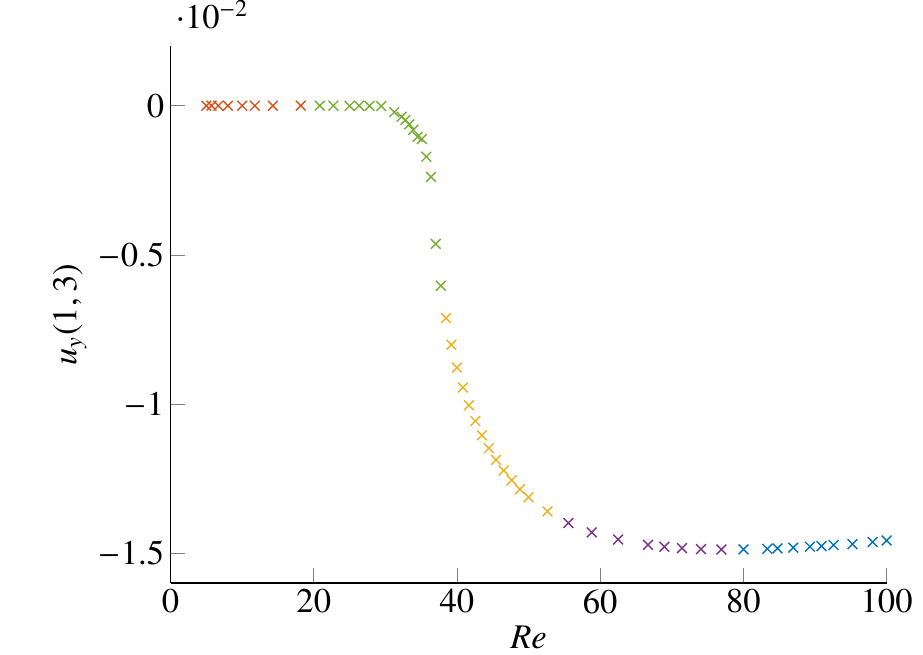}
\includegraphics[scale=0.45]{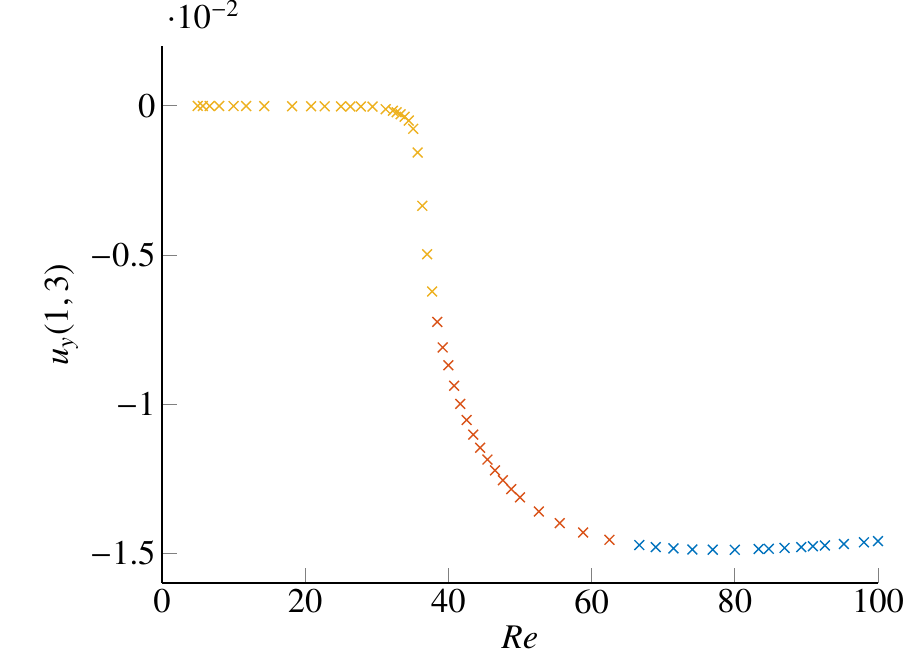}
\includegraphics[scale=0.45]{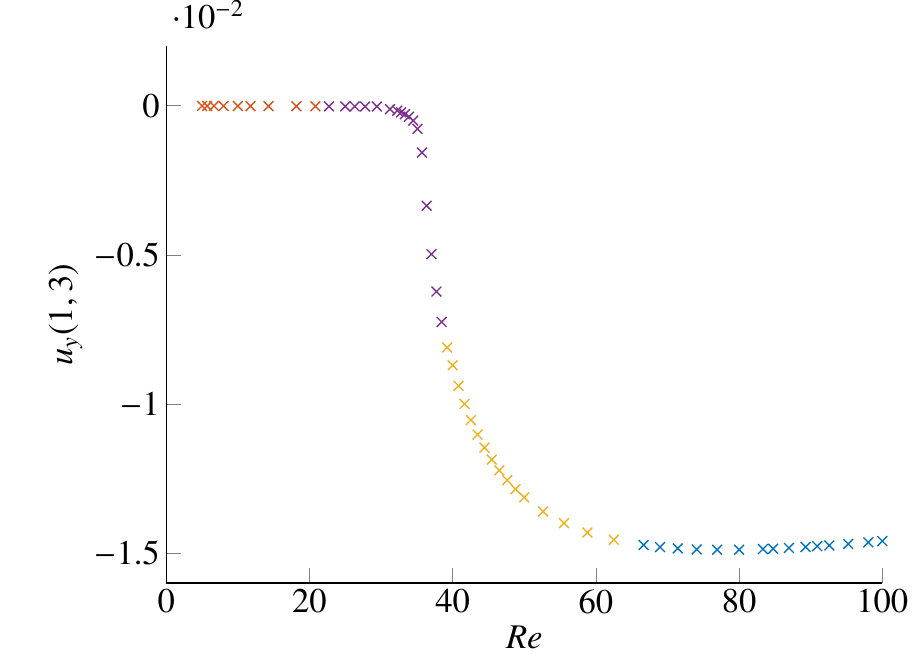}
\includegraphics[scale=0.45]{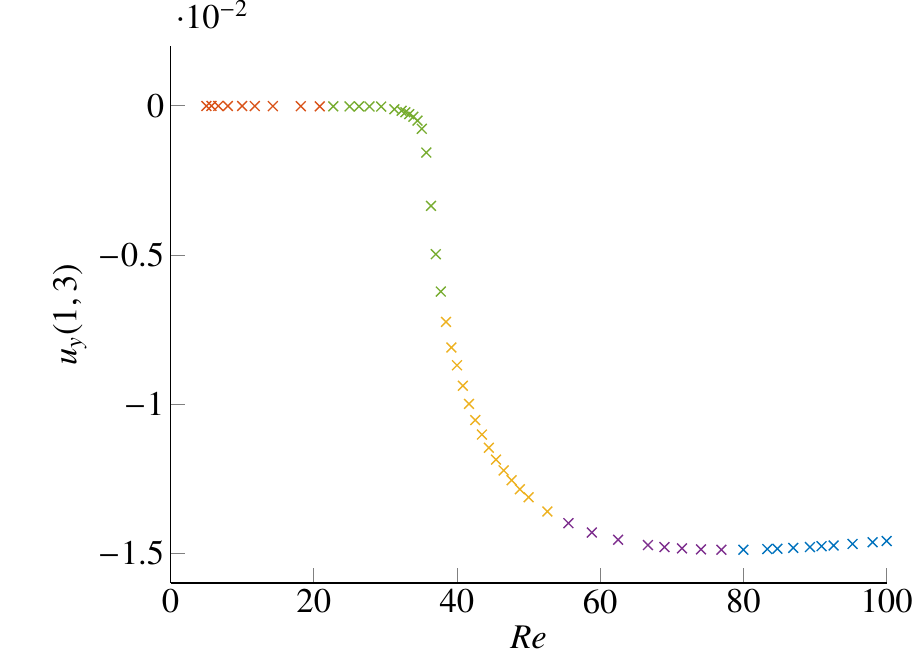}
\end{center}
\caption{Left to right:  k-means clusterings of snapshots for $3$, $4$, and $5$ clusters correspond to the use of the SEM (top row) and FEM (bottom row) methods.}
\label{fig:cluster_comparison}
\end{figure}

\vskip5pt
\noindent\textbf{\em Clustering.}
We apply the k-means method to cluster the snapshots produced by the two spatial discretization approaches. Results are shown in Fig.~\ref{fig:cluster_comparison} for $3$, $4$, and $5$ clusters. The k-means method decides on how to cluster the snapshots by minimizing the total clustering variance. Fig. \ref{fig:Channel_CVT} hints that two discretization methods will result in very similar clusterings; this is confirmed in Fig. \ref{fig:cluster_comparison} for which one observes very close agreement between the clusterings obtained using the two methods.

\subsubsection{Results for the online stage: the use of the ROM bases with local bases}
Not surprisingly given the close agreement between the snapshot results obtained using the two spatial discretization methods, their application to the use of ROMs also result in very close agreements. As described in \S\ref{sec:bas}, we use the clustering of the snapshots to produce a POD basis for each cluster. Having those bases in hand, we then use them to determine ROM solutions, i.e., solutions of \eqref{abs-dis2} using local POD bases. Of course, we do so for parameter values, i.e., Reynolds numbers in the current setting, not used in the determination of snapshots. Both SEM and FEM ROM solutions were obtained using an Oseen fixed-point iteration.


%
%
%

Given a new parameter value, one must decide which local basis to use. In \S\ref{sec:bif_detection}, two criteria where discussed for this purpose: using the mean of the clusters or using a midrange/radius approach.  
The difference in the approaches is that the parameter centroid takes the distribution of snapshots in parameter space into account whereas the mid-range uses only information about the extent of the cluster in parameter space. 
Local ROMs are computed by POD using the snapshots belonging to their respective clusters. To recover the bifurcation diagram, the parameter mean and the midrange/radius criteria are used to decide which local basis to use for a given parameter value of interest. 

\begin{figure}[h!]
\begin{center}
\includegraphics[scale=.45]{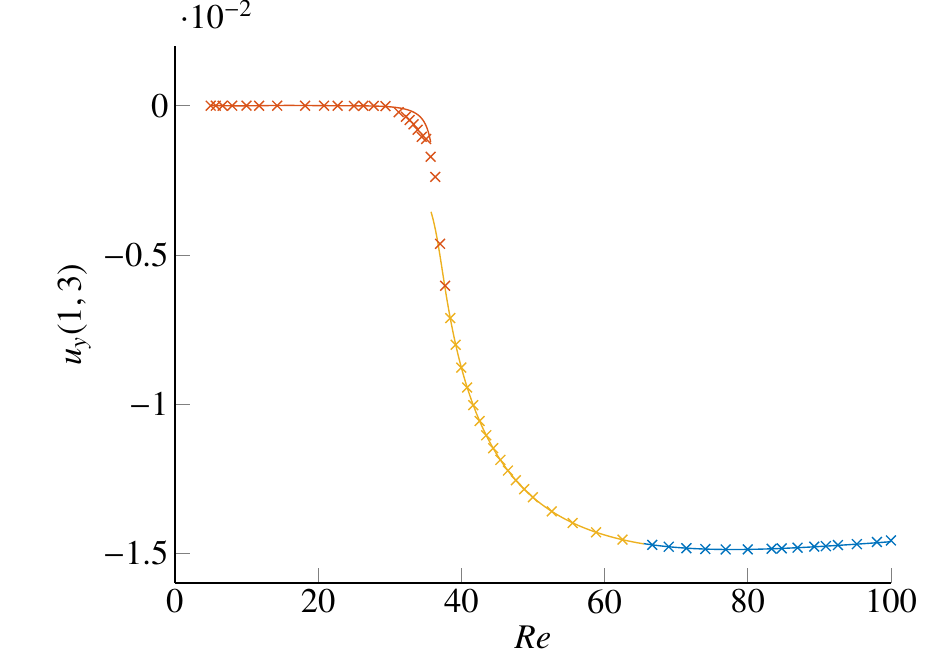}
\includegraphics[scale=.45]{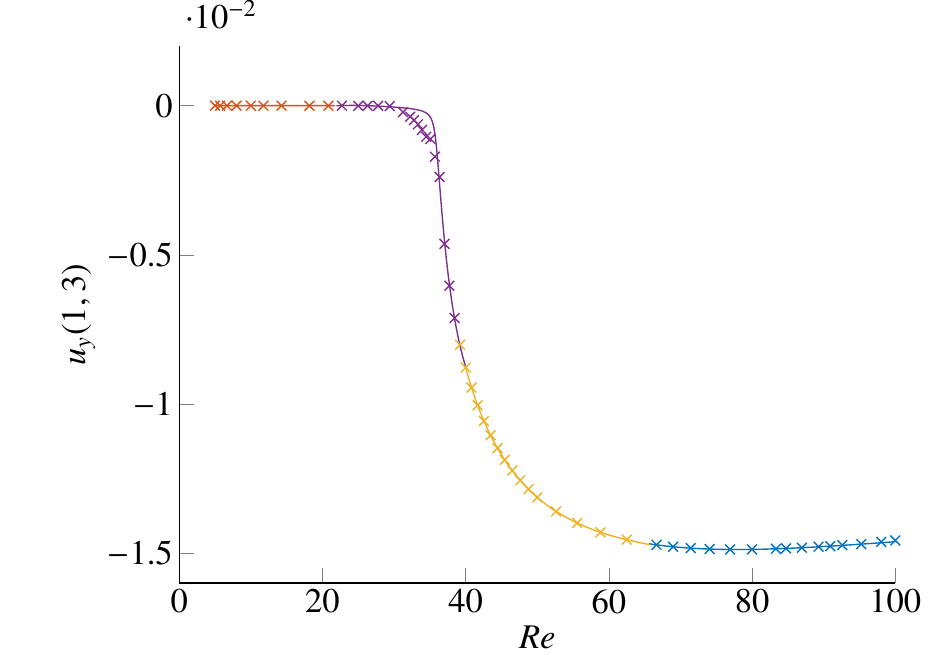}
\includegraphics[scale=.45]{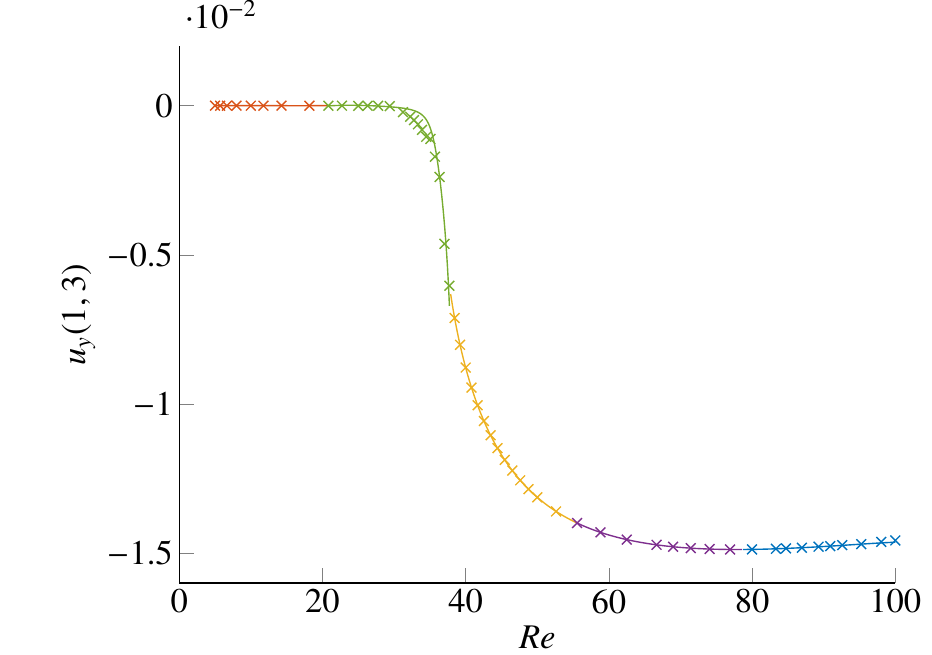}
\includegraphics[scale=.45]{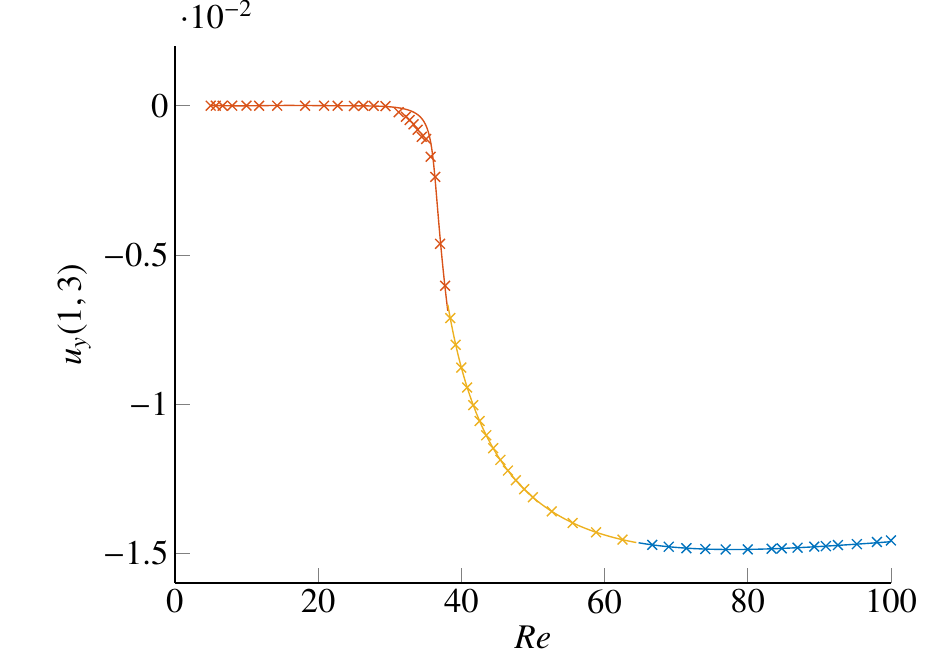}
\includegraphics[scale=.45]{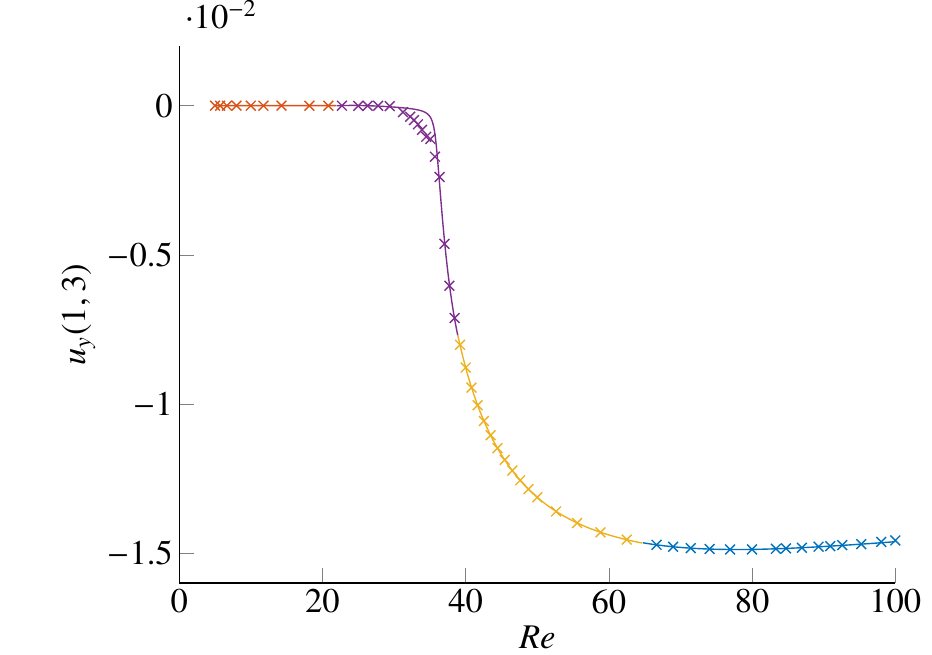}
\includegraphics[scale=.45]{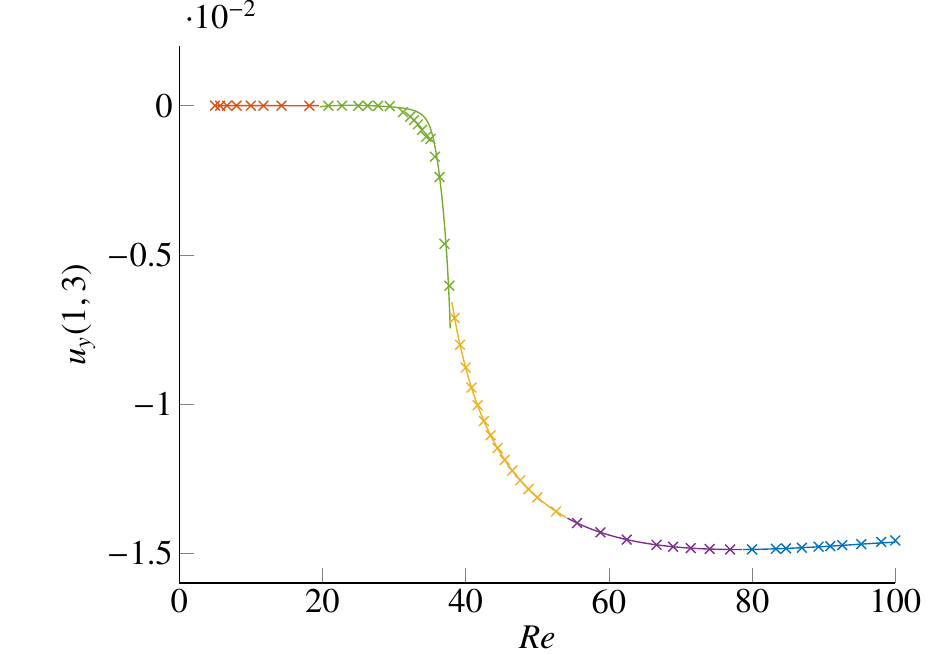}
\caption{For the SEM spatial discretization and for 3, 4, and 5 clusters (left to right), bifurcation diagrams and assignments of parameters to local ROMs using the distance to parameter mean criterion (top row) and using the distance to parameter mid-range and cluster radius criterion (second row). The crosses correspond to full-order solutions for parameter values not used to generate the snapshots whereas the solid curves correspond to ROM solutions.} \label{fig:Channel_3_4_5c}
\end{center}
\end{figure}

\begin{figure}[h!]
\begin{center}
\includegraphics[scale=.45]{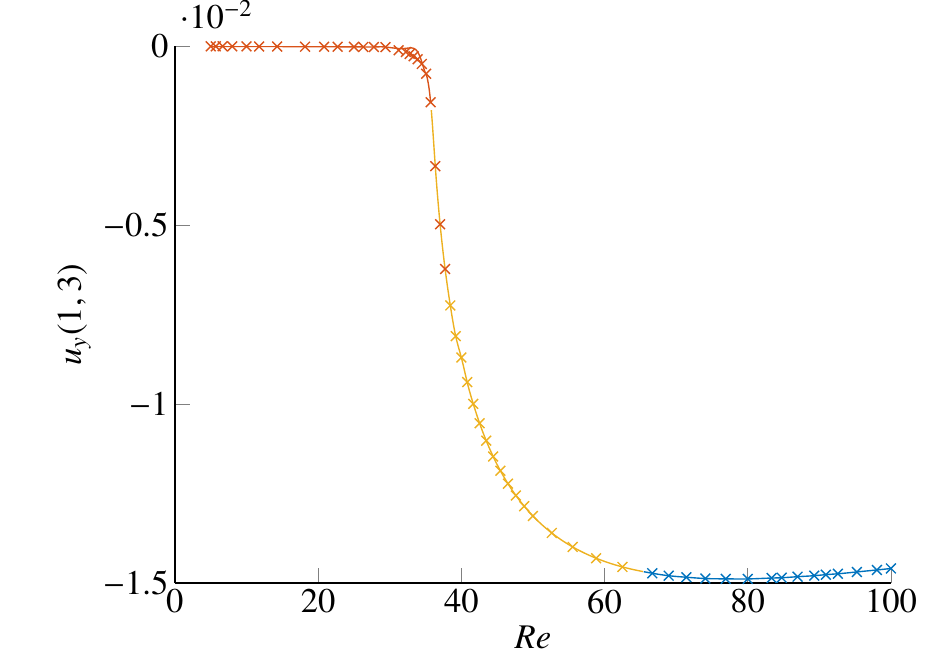}
\includegraphics[scale=.45]{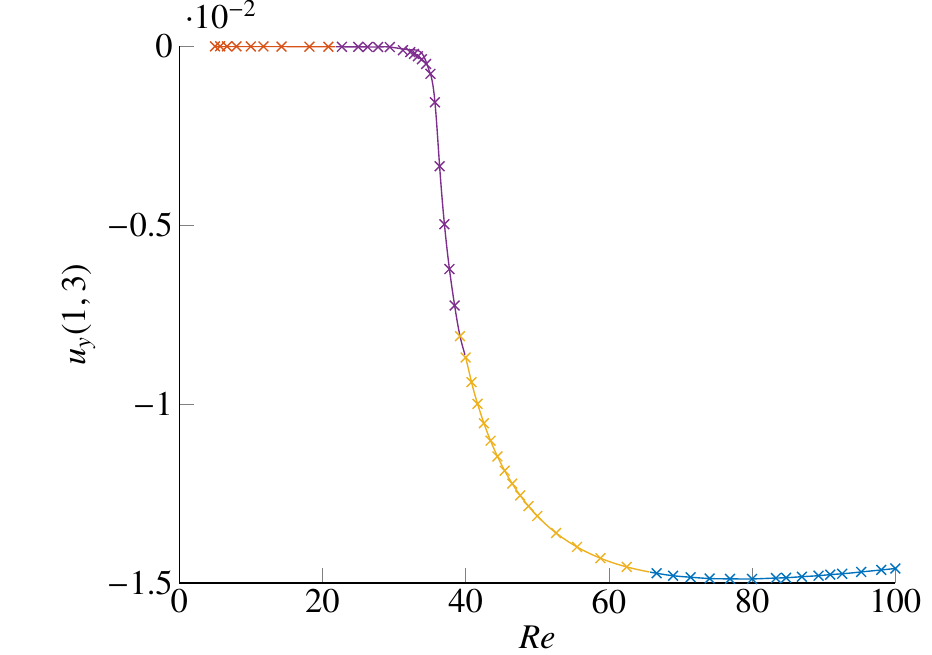}
\includegraphics[scale=.45]{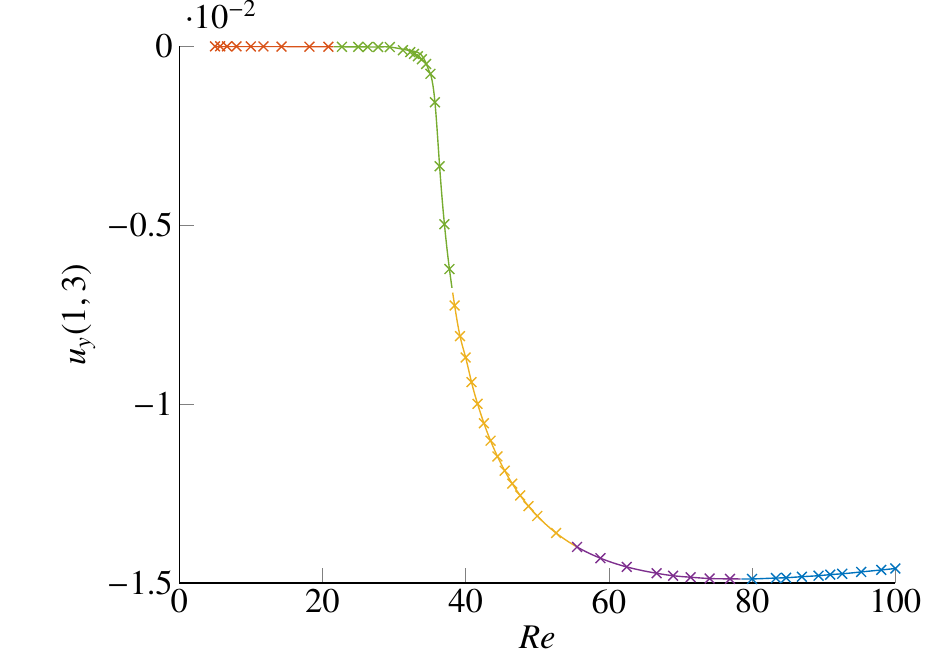}
\includegraphics[scale=.45]{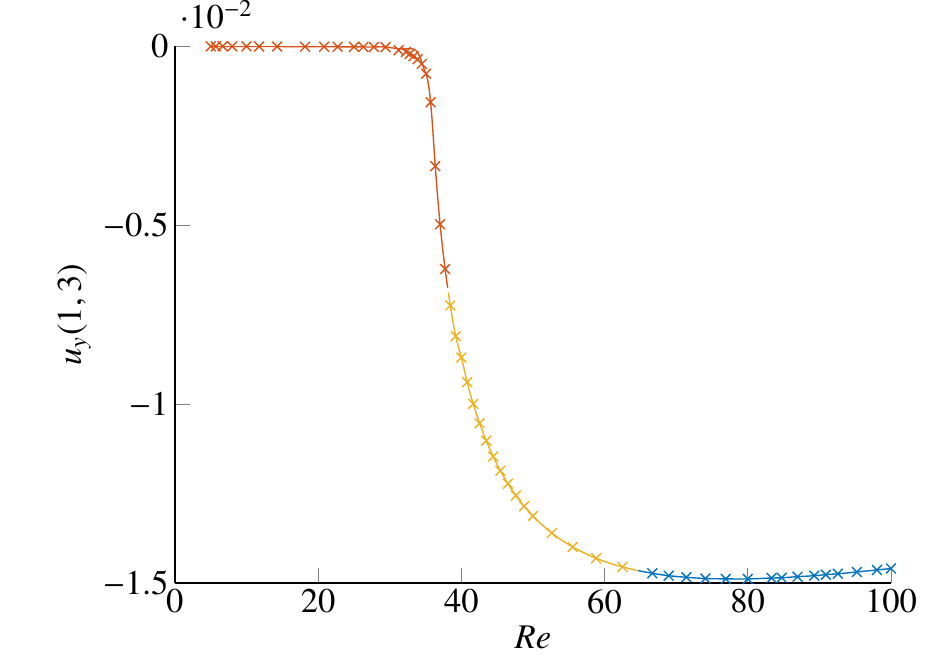}
\includegraphics[scale=.45]{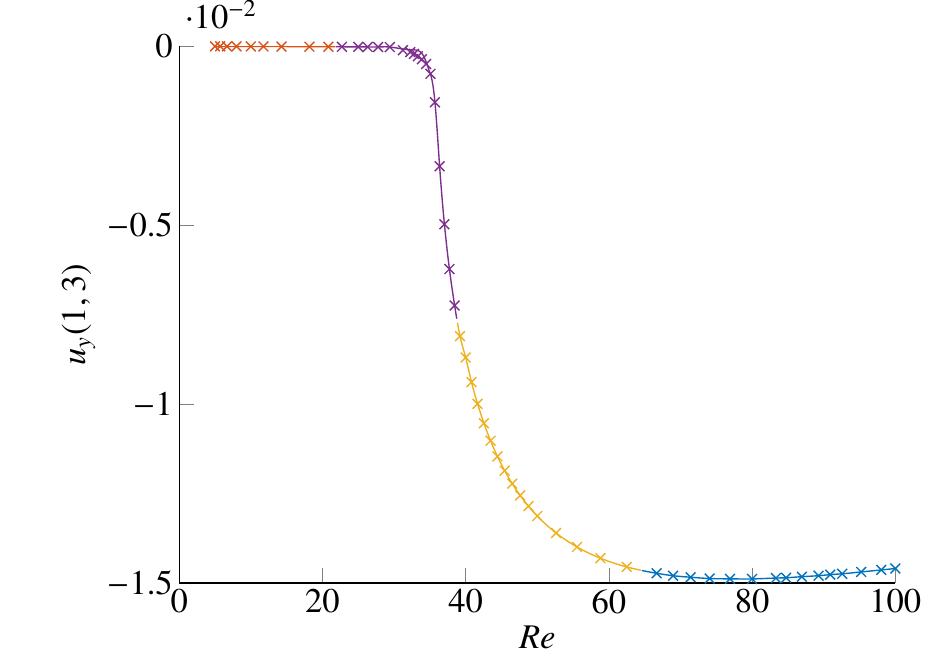}
\includegraphics[scale=.45]{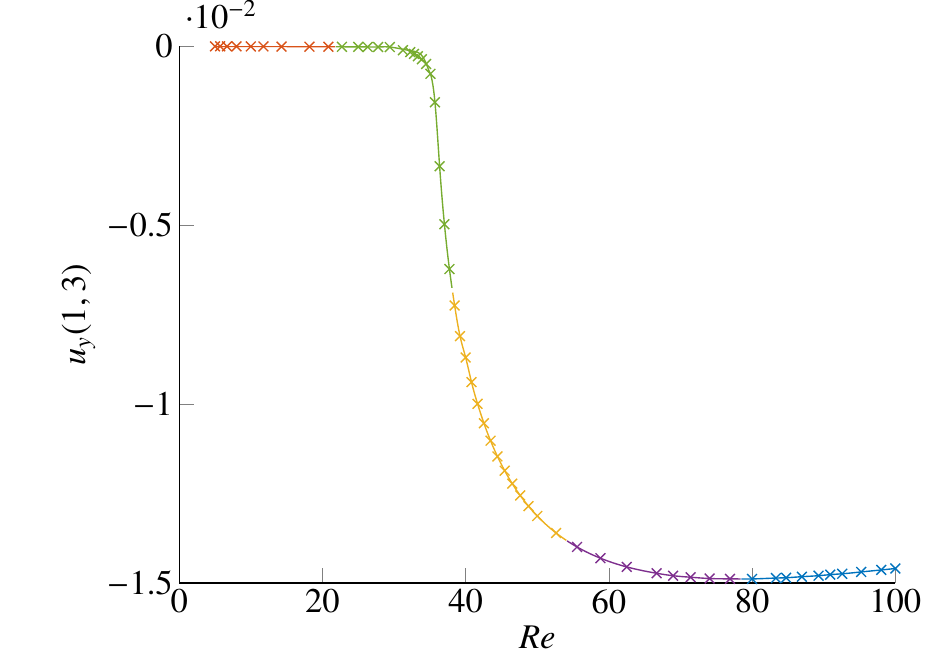}
\caption{The same information as in Figure \ref{fig:Channel_3_4_5c} for the FEM spatial discretization.} 
  \label{fig:FEM_Channel_3_4_5F}
\end{center}
\end{figure}

In Figure~\ref{fig:Channel_3_4_5c}, bifurcation diagrams and assignments of parameters to local ROM bases are given for 3, 4, and 5 clusters with respect to SEM snapshots (left to right).  The rows correspond to the use of the distance to parameter mean criteria (top row) and using the distance to parameter mid-range and cluster radius criteria (bottom row). The crosses correspond to full-order solutions for parameter values not used to generate the snapshots; they are color coded to indicate which cluster they belong to. The solid curves correspond to to ROM solutions using many more new parameter values with the color coding corresponding to which cluster, and therefore which local ROM basis, a new parameter is assigned to. We observe that the midrange/radius criterion assigns new parameters to the correct local basis. On the other hand, mis-assignments of local bases occurs when the parameter mean criteria is used. The same information for local ROMs relative to FEM discretization are given in Figure \ref{fig:FEM_Channel_3_4_5F}. We note that the results are in a excellent agreement with those in Figure  \ref{fig:Channel_3_4_5c}.

\vskip5pt
\noindent\textbf{\em Errors with respect to full-order models and comparisons to global ROM models.} The goal of reduced order modeling is to obtain accurate approximations of full-order model solutions while incurring significantly reduced compute costs. In addition, because we propose a localized ROM basis approach based on clustering, comparison of the accuracy and compute times for using local ROMs should also be made with respect to a global ROM. Quantitive comparisons are made with respect to the relative error
\begin{equation}\label{relerr}
\mbox{\em error} = \frac{\|u_{\mbox{\footnotesize\em full-order}} - u_{\mbox{\footnotesize\em rom}}\|_{L^2(\Omega)}}{\|u_{\mbox{\footnotesize\em full-order}}\|_{L^2(\Omega)}},  
\end{equation}
where $u_{\mbox{\footnotesize\em rom}}$ could correspond to local or global ROMs. An investigation of ROM stability and error certification with respect to the Navier-Stokes equations can be found in \cite{Yano:2013}.

We note that the offline computational cost is dominated by constructing the discrete trilinear form appearing in the weak form of the Navier-Stokes equations and used to project the nonlinearity; that cost grows cubically in the reduced basis size $L$. The online computation costs are negligible by comparison, because no operations involving the full-order model dimension are necessary.

Table \ref{jet_3clust_accu} provides, for 3, 5, and 7 clusters, the relative error between the full-order solutions and the local and global ROM solutions as well as the size of the bases. Ten values of the Reynolds number are used to determine solutions of the full-order model, the local ROM, and the two global ROMs. Although lying within the Reynolds number range of interest, the Reynolds number values chosen are all different from those used to generate the snapshots used to construct the ROMs. One global ROM, referred to as Global-1 in the table, uses a basis of dimension equal to the sum the dimensions of all the local ROM sizes. Whereas this global ROM has a larger compute time than the local ROM approach, it can serve as a baseline reference for assessing the approximation quality that local ROMs provide. The second global ROM, referred to as Global-2 in the table, uses a basis of dimension equal to the maximum dimension of the local bases used so that it has a similar offline computational time as that of the local ROM approach. For the local bases, results are given for both the use of the cluster mean and midrange/radius criteria when selecting what local basis to use. Values marked with (n.c.) correspond to non-converged cases; the iteration is terminated after $3000$ iterates, whereas convergence typically happens after about $100$ iterations. In this case, we report the relative error value for last iterate computed. Also, the arithmetic mean and maximum of relative errors (non-converged values are included) over the 10 realizations are also given in the tables. 

\begin{table}[h!]\caption{For the jet problem: cluster accuracy and basis size comparisons.} \label{jet_3clust_accu}
\begin{tiny}
\begin{center}
\begin{tabular}{l*{6}{c}c}
\multicolumn{6}{l}{\underline{\bf $3$ clusters}}\\[.5ex]
\multicolumn{1}{c}{ROM}    & Re=6 & Re=16.333 & Re=26.666 & Re=37 & Re=47.333 & Re=56.666 \\
\hline
Global-1    & 0.00054 & 0.00052 & 0.00052 & 0.0038 & 0.00079 & 0.945 (n.c.) \\
Global-2    & 0.0027 & 0.0016 & 0.00068 & 0.0043 & 0.00091 & 0.00099 \\
Local    & 0.00081 & 0.00171 & 0.00121 & 0.011/0.021 & 0.00081 & 0.0020 \\[1ex]
\multicolumn{1}{c}{\em ROM}     & Re=68 & Re=78.333 & Re=88.666 & Re=99 & mean & max & basis size\\
\hline
Global-1    & 1.02 (n.c.) & 1.00 (n.c.) & 0.996 (n.c.) & 0.995 (n.c.) & 0.50 & 1.02 & 39 \\
Global-2   & 0.0013 & 0.14 (n.c.) & 0.26 (n.c.) & 0.67 (n.c.) & 0.11  & 0.67 & 21 \\
Local ROM   & 0.0013 & 0.0084 & 0.0050 & 0.0167 & 0.0064 & 0.0167  & 11/7/21 \\[2ex]
\multicolumn{6}{l}{\underline{\bf $5$ clusters}; *=10/8/5/14/17}
\\[.5ex]
\multicolumn{1}{c}{ROM}    & Re=6 & Re=16.333 & Re=26.666 & Re=37 & Re=47.333 & Re=56.666 \\
\hline
Global-1 & 0.00054 & 0.00052 & 0.00052 & 0.0037 & 0.00079 & 0.00078  \\
Global-2 & 0.0043 & 0.02 & 0.0082 & 0.0064 & 0.0017 & 0.0019 \\
Local    & 0.00054 & 0.00052 & 0.00060 & 0.0285 & 0.015 & 0.00078\\[1ex]
\multicolumn{1}{c}{ROM}    & Re=68 & Re=78.333 & Re=88.666 & Re=99 & mean & max & basis size\\
\hline
Global-1   & 0.1415 (n.c)  & 0.240 (n.c) & 0.611 (n.c) & 0.829 (n.c) & 0.18 & 0.829 & 54 \\
Global-2   & 0.0022 & 0.178 (n.c) & 1.64 (n.c) & 1.08 (n.c)    & 0.29 & 1.64 & 17 \\
Local      & 0.0013 & 0.0084 & 0.005 & 0.045 & 0.011 &  0.045 & *\\[2ex] 
\multicolumn{6}{l}{\underline{\bf $7$ clusters}; **=8/6/4/4/8/6/16}\\[.5ex]
\multicolumn{1}{c}{ROM}    & Re=6 & Re=16.333 & Re=26.666 & Re=37 & Re=47.333 & Re=56.666 \\
\hline
Global-1 & 0.00054 & 0.00052 & 0.00052 & 0.0037 & 0.00079 & 0.00079  \\
Global-2 & 0.016 & 0.0077 & 0.0052 & 0.018 & 0.0039 & 0.0068 \\
Local    & 0.00054 & 0.00052 & 0.00057 &  0.0036 & 0.0017 & 0.0014\\[1ex]
\multicolumn{1}{c}{ROM}   & Re=68 & Re=78.333 & Re=88.666 & Re=99 & mean & max & basis size\\
\hline
Global-1   & 0.152 (n.c)  & 0.23 (n.c) & 0.619 (n.c) & 0.712 (n.c) & 0.17 & 0.712 & 52 \\
Global-2   & 0.41 (n.c.) & 0.48 (n.c) & 7.8 (n.c) & 1.7 (n.c)    & 1.04 & 7.8 & 16 \\
Local      & 0.0013 & 0.0083 & 0.0050 & 0.017 & 0.004 & 0.017  & ** 
\end{tabular}
\end{center}
\end{tiny}
\end{table}

For the 3-cluster case, the local POD basis sizes are $11$, $7$, and $21$, respectively, with $11$ corresponding the smaller values of Re and 21 to the larger values. The two global basis dimensions are thus $39$ and $21$. The switch between local ROMs occurs at Re = 35 and Re = 65 for the cluster mean criterion and at Re = 38 and Re = 64 for the midrange/radius criterion. Thus, two error values are shown for Re = 37. For the 5-cluster case, the switches between local ROMs occurs at Re = 21, 38, 56, and 78 for the parameter cluster mean criterion and Re = 19.5, 38.1, 54.1, and Re = 78.5 for the cluster midrange/radius criterion. For the 7-cluster case, the switches between local ROMs occurs at Re = 20.2, 34.9, 42.7, 51.7, 65.2, and 82.1 for the parameter cluster mean criterion and  at Re = 19.5, 36.7, 42.1, 51.3, 64.6, and 81.7 for the cluster midrange/radius criterion.

The tabulated results suggest that the use of local reduced bases provide a more efficient and usually more accurate ROM approach when compared to the use of global reduced bases. For the one case in which the two different local basis selection criteria differed (see Table \ref{jet_3clust_accu}), the parameter cluster midrange/radius criterion resulted in a smaller error.


\subsection{Example 2: Rayleigh-B\'{e}nard cavity flow}\label{sec:cav}

As a second example, we consider Rayleigh-B\'{e}nard (RB) cavity flow, 
introduced in~\cite{Roux:GAMM} and widely studied since then. This test case features a rich bifurcating behavior despite the very simple geometrical setting. The domain $\Omega$ is a rectangular cavity with height $1$ and length $A$. We choose the bottom left corner of the cavity as the origin of the coordinate system. The vertical cavity walls are maintained at constant temperatures $T_0$ (left wall) and $T_0 + \Delta T$ (right wall) 
with $\Delta T > 0$, whereas the horizontal walls are thermally insulated. A convective flow is induced by the horizontal component of the temperature gradient. The fluid is assumed to have density $\rho$ and viscosity $\nu$ at the reference temperature $T_0$ and to satisfy homogeneous boundary condition on the cavity walls. 

The strong formulation of the RB cavity problem, as stated in~\cite{Roux:GAMM}, 
is given by \eqref{NS-1} with $\f = (0, g \beta (T - T_0) (x/A) \e_y)^T = (0, g \beta \Delta T (x/A) \e_y)^T$. Here, $g$ denotes the magnitude of the gravitational acceleration, $\beta$ is the coefficient of thermal expansion, 
$x$ is the horizontal coordinate, and $\e_y$ is the unit vector directed along the vertical axis. The Grashof number
\begin{equation}\label{eq:gr}
\text{Gr} = \frac{g \beta \Delta T}{A \nu^2}
\end{equation}
can be used to characterize the flow regime. This definition already takes into account the fact that the height of the cavity is $1$. The Grashof number can be thought of as the ratio of buoyancy forces to viscous forces. For large Grashof numbers buoyancy forces are dominant over viscous forces and vice versa. Note that with \eqref{eq:gr} we can write  ${\bm f} =(0, \text{Gr}\nu^2 x)^T$.
The Prandtl number for this problem is zero. Most of the numerous numerical studies of this problem consider $A = 4$ as was done in \cite{Roux:GAMM}.

Depending on the governing parameters and the initial state, three branches of steady state solutions exist with one, two, or three primary circulations, also referred to as rolls. If the initial state is fluid at rest, for sufficiently small value of $\text{Gr}$, a steady single roll flow is observed.  For a value of $\text{Gr}$ greater than a first critical value, the steady-state velocity field features two rolls. The two-roll flow is a stable configuration until $\text{Gr}$ exceed a second critical value, at which the steady-state flow has three rolls. Snapshots of the flow fields for three stable scenarios are illustrated in \S\ref{sec331}.  The passage of $\text{Gr}$ from one subregion of the bifurcation diagram to another causes a {\em discontinuous change} in that the number of rolls changes. 
Note that all steady-state flows are centrally symmetric, i.e., with respect to rotation through $180^\circ$ about the vertical centerline of the cavity, and also note that all rolls rotate in a counter-clockwise direction.

The model is defined on a rectangular domain $\Omega = [0, 4] \times [0, 1]$, i.e., the aspect ratio is set to $A = 4$. The Grashof number $Gr$ is considered as a parameter ranging over the interval  $[10^{3}, 10^{5}]$; the viscosity $\nu$ is set to one.

\begin{figure}[h!]
\begin{center}
 \includegraphics[scale=.15]{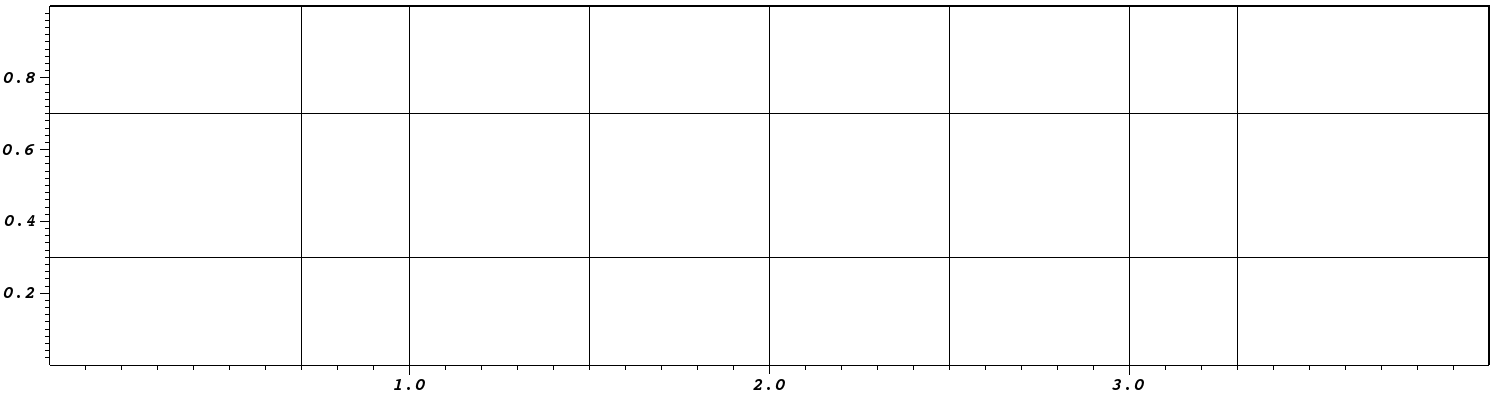} 
 \caption{The $24$ spectral elements used in the simulations, resulting in $18963$ degrees of freedom.}
 \label{Hess:domain_cavity}
\end{center}
\end{figure}

Because of the close agreement between the results obtained using the SEM and FEM spatial discretizations for the jet problem, for the sake of brevity, for the cavity problem we only report on numerical results obtained for the SEM discretization. The cavity domain is subdivided into $24$ quadrilateral elements as shown in Fig.~\ref{Hess:domain_cavity} and modal Legendre ansatz functions  of order $16$ are used in every element and solution component. 
This leads to $6321$ global degrees of freedom in each component, i.e., horizontal velocity, vertical velocity, and pressure for the time-dependent simulation.

\noindent\textbf{Timings.}
The full-order solutions are obtained using time-marching scheme with the C++ soltware Nektar++. On the other hand, the 
reduced-order solutions are obtained using a fixed point iteration scheme implemented in python. As a result,  a direct comparison of timings is not possible. 
However, comparing a single step of the fixed-point scheme using for the  full-order problem with the online reduced-order problem, compute times reduce
from $100$s to $0.1$s.

\subsubsection{Results for the offline stage: the construction of local ROM bases}\label{sec331}
In this section, we provide numerical results illustrating the different aspects of local ROM construction process for the cavity problem.

\begin{figure}[h!]
\begin{center}
\includegraphics[width=2.25in]{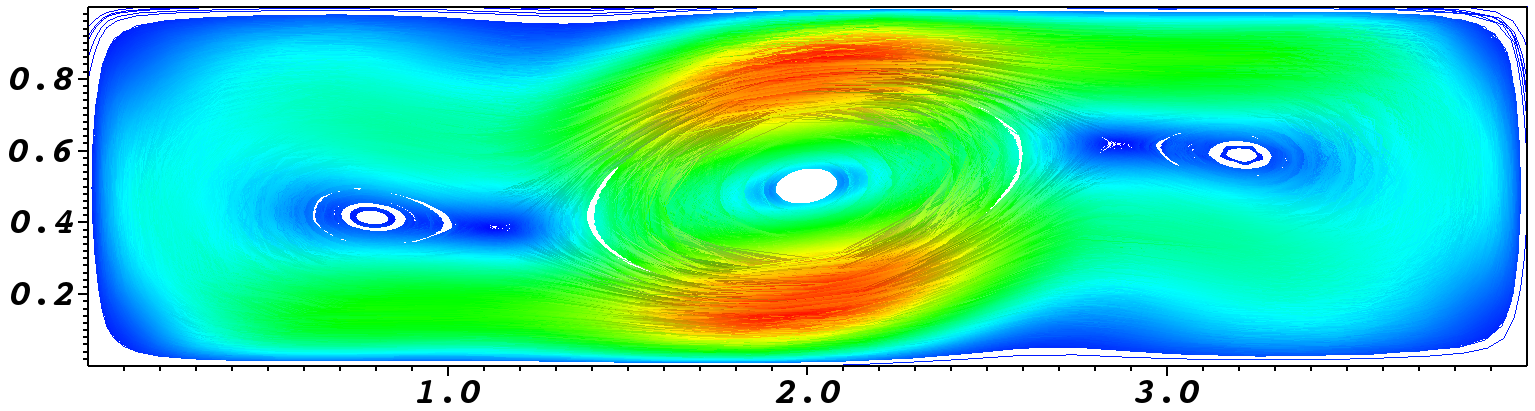}
\includegraphics[width=.7in]{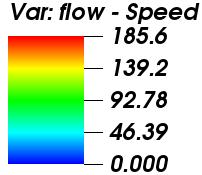}
\vskip5pt
\includegraphics[width=2.25in]{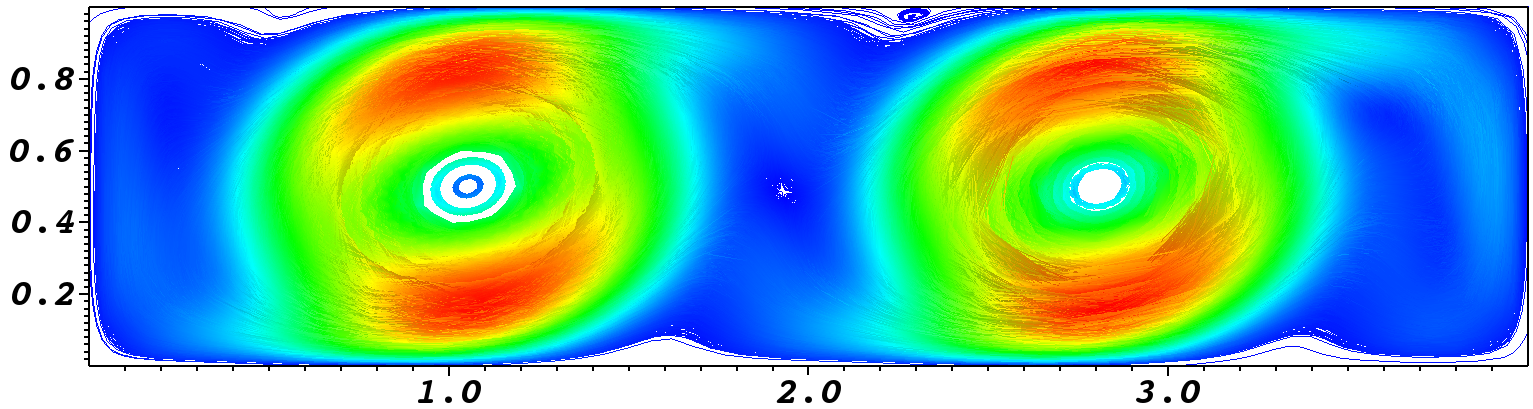}
\includegraphics[width=.7in]{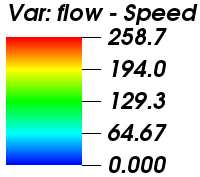}
\vskip5pt
\includegraphics[width=2.25in]{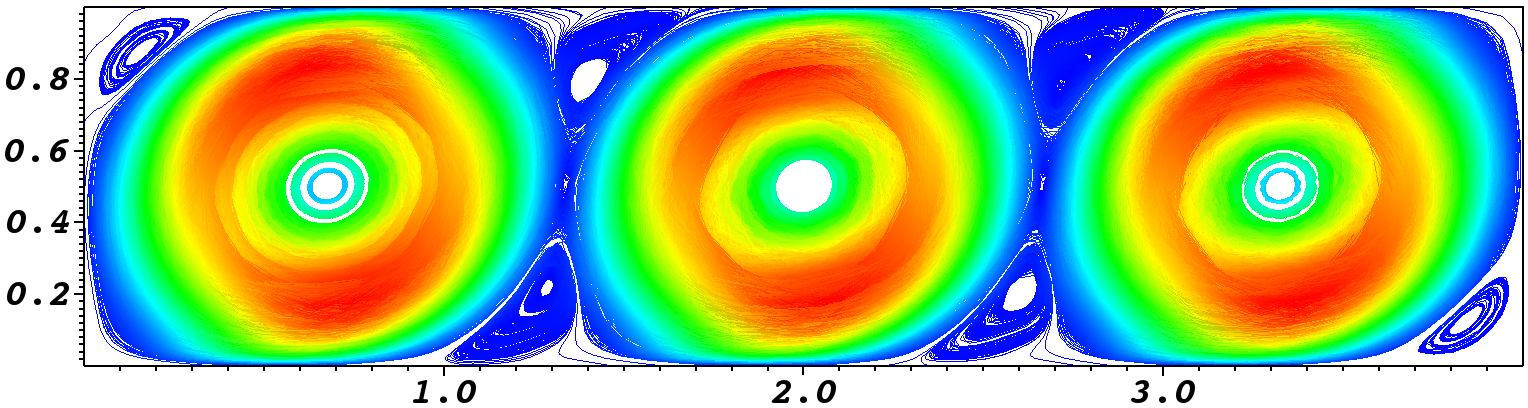}
\includegraphics[width=.7in]{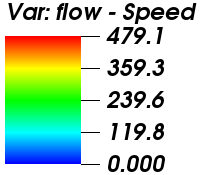}
\caption{Steady state solutions of the cavity flow problem for Grashof numbers $20 \cdot 10^{3}$, $40 \cdot 10^{3}$, and $98 \cdot 10^{3}$ (top to bottom).} \label{fig:cavity_20_sol}
\end{center}
\end{figure}

\vskip5pt
\noindent\textbf{\em Snapshots.} Fig.~\ref{fig:cavity_20_sol} shows representative 
steady state solutions of this model for Grashof number Gr $=20 \cdot 10^{3}$, $40 \cdot 10^{3}$, and $98 \cdot 10^{3}$, respectively. For {\text Gr} up to $30 \cdot 10^{3}$, steady state solutions feature a single roll, whereas two rolls develop after that. For {\text Gr} in the interval between $92 \cdot 10^{3}$ and $98 \cdot 10^{3}$, three rolls are observed.

There also exist time-periodic solutions. For example, for Gr  $=40 \cdot 10^{3}$, in addition to the steady state solution with two rolls, a time-periodic solution with one roll exists.  Depending on the initial condition, either solution
can be obtained. For ${\text Gr}=100 \cdot 10^{3}$ a time-periodic solution with three rolls also exists.

To compute a steady state, the initial condition is either given as zero or a steady state at a close parameter location. The time-step size is set to $10^{-7}$ and about $10^8$ time steps are required to reach a steady state. Steady state solutions with one roll were found for Gr $\in[ 10 \cdot 10^{3}, 22 \cdot 10^{3} ]$, two rolls for Gr $\in[ 30 \cdot 10^{3}, 90 \cdot 10^{3} ]$, and three rolls for Gr $\in[ 92 \cdot 10^{3}, 98 \cdot 10^{3} ]$. 

\vskip5pt

\begin{figure}[h!]
  \begin{center} 
\includegraphics[width=.45\textwidth]{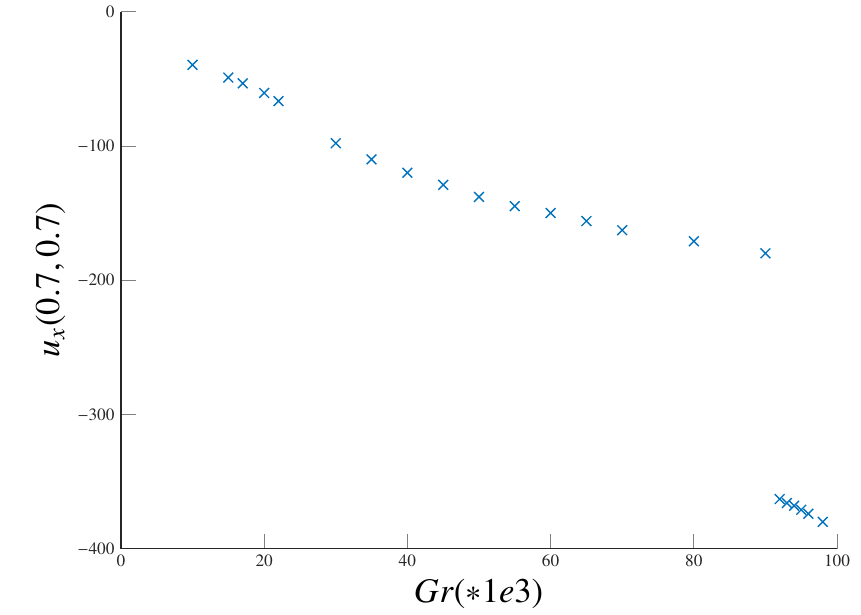} 
    \end{center}   
\caption{Bifurcation diagram for the cavity problem corresponding to the horizontal velocity component at $(0.7, 0.7)$.} \label{fig:cavity_bif_dia}   
\end{figure}

\noindent\textbf{\em Bifurcation diagram.} Using full-order solutions, e.g., snapshots, we have the bifurcation diagram corresponding to the point value of the horizontal component of the velocity at $(0.7, 0.7)$  shown in Fig.~\ref{fig:cavity_bif_dia}. The two bifurcations are clearly evident, although the one at lower Gr number is ``less severe'' than the one at higher Gr number.

\begin{figure}[h!]
\begin{center}
\includegraphics[width=0.45\textwidth]{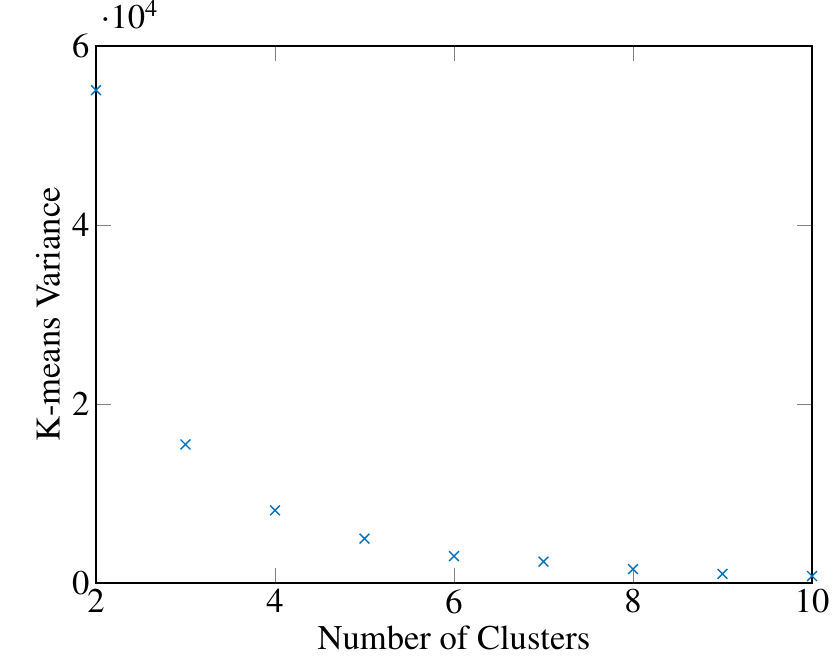}  
\end{center}
\caption{k-means variance versus number of clusters for the cavity problem.} \label{fig:cavity_CVT_silhouette}  
\end{figure}

\begin{figure}[h!]
\begin{center}
\includegraphics[scale=0.4]{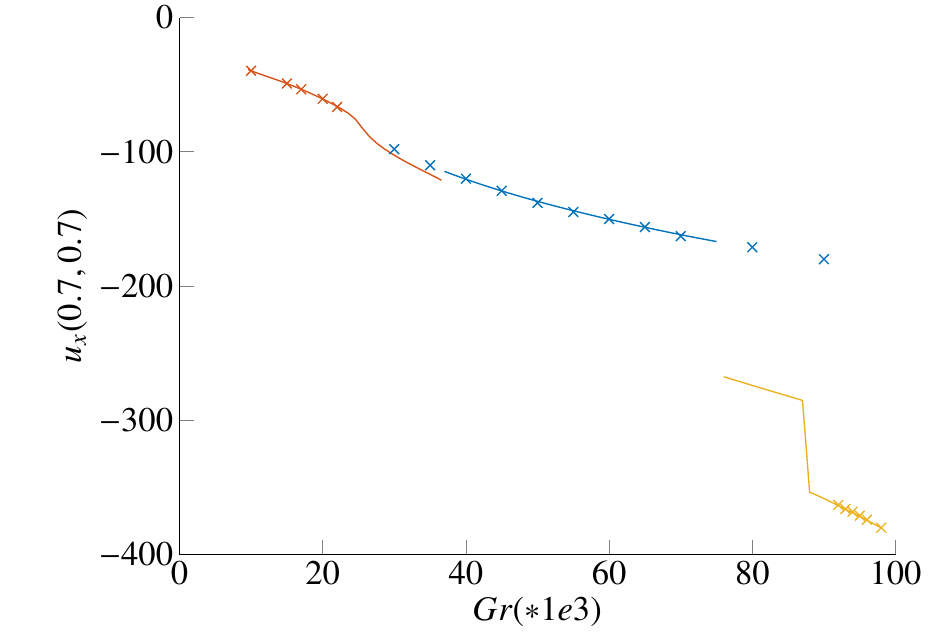}
\includegraphics[scale=0.4]{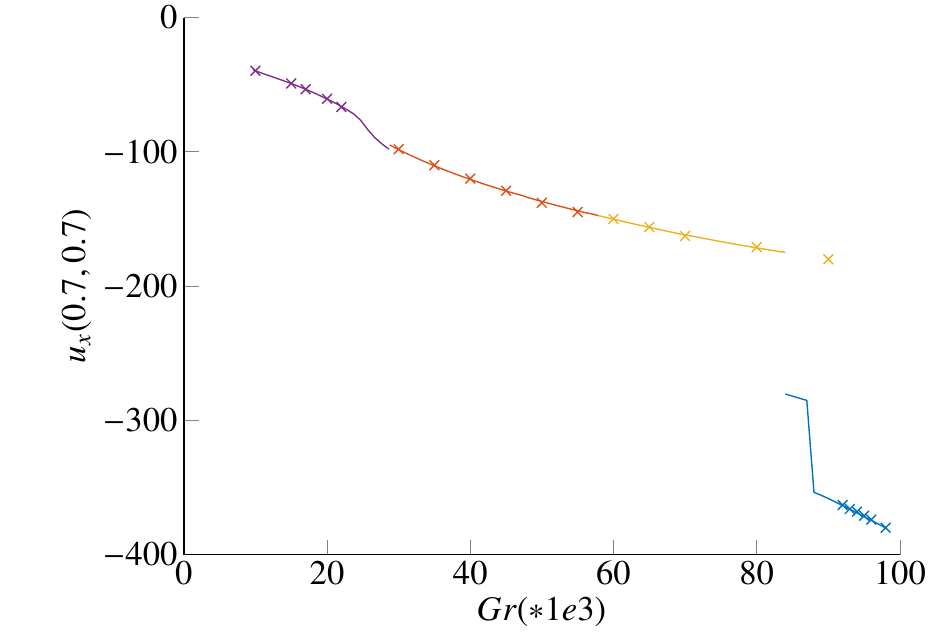}
\includegraphics[scale=0.4]{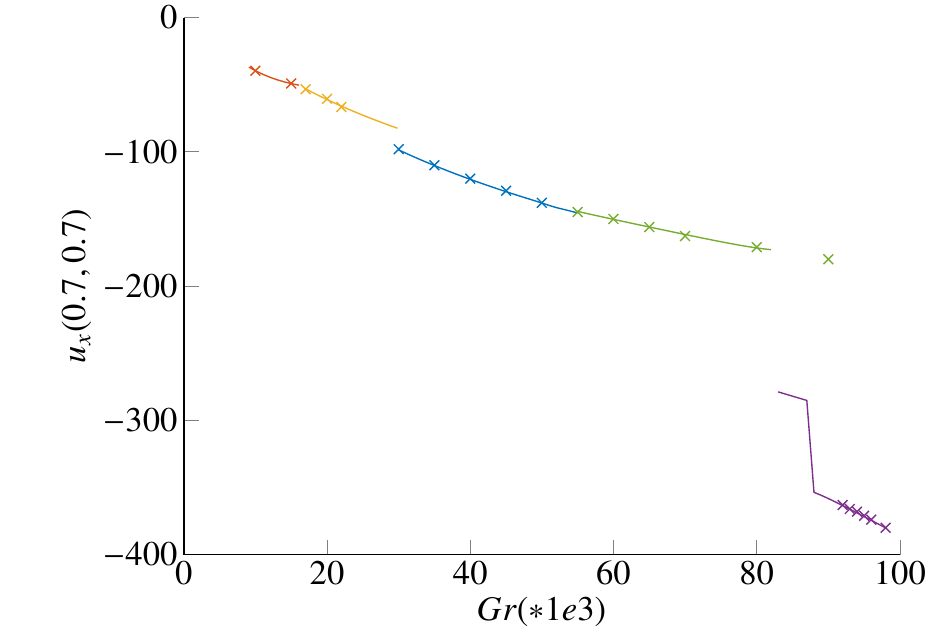}
\includegraphics[scale=0.4]{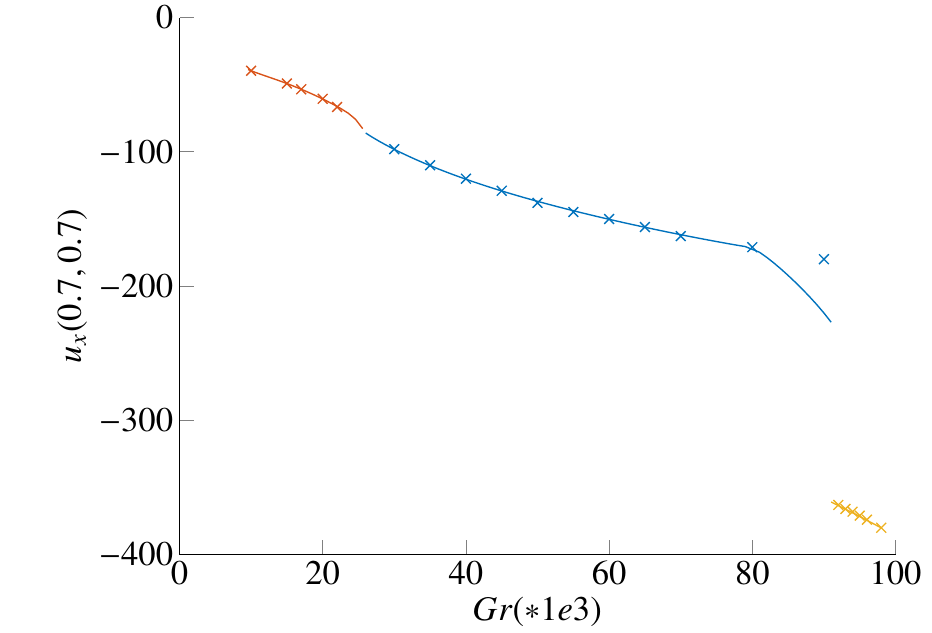}
\includegraphics[scale=0.4]{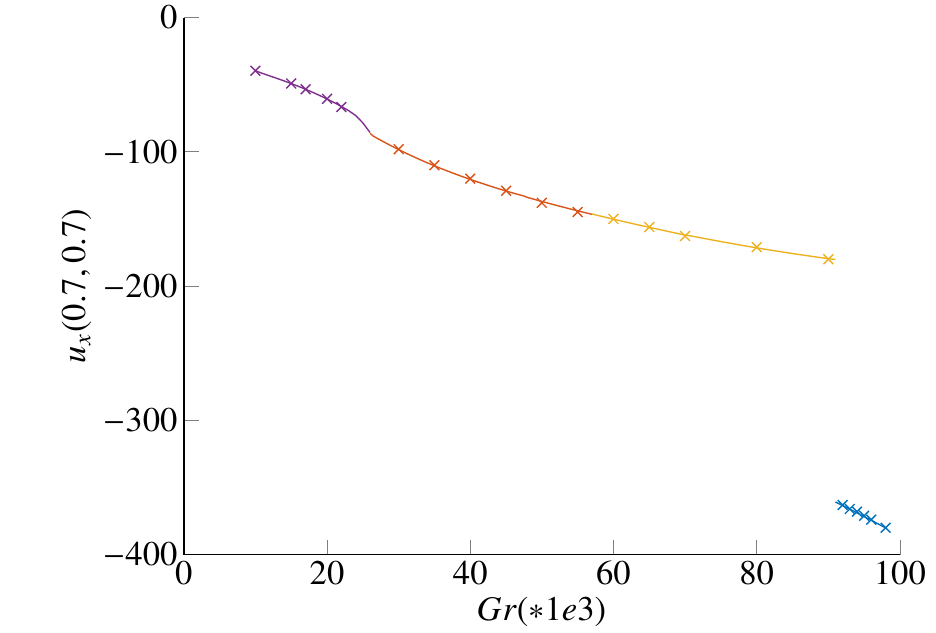}
\includegraphics[scale=0.4]{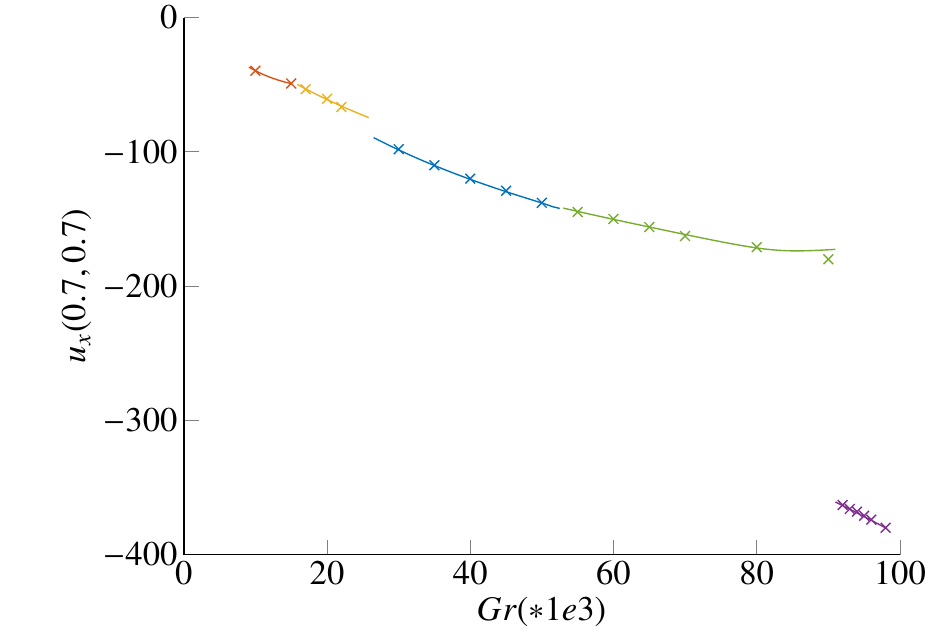}
\caption{Recovered bifurcation diagram with 3 (left), 4 (middle), and 5 (right) clusters using the distance to parameter mean criterion (first row) and using the distance to parameter mid-range and cluster radius criterion (second row).}
\label{fig:cavity_3_4_5c}
\end{center}
\end{figure}

\vskip5pt
\noindent\textbf{\em Clustering.}
In Fig.~\ref{fig:cavity_CVT_silhouette}, the k-means variance is plotted vs. the number of clusters in the range $2$ to $10$. We observe a faster reduction in the variance compared to what we saw for the jet problem; compare with Fig. \ref{fig:Channel_CVT}. The clustering algorithm separates the one, two, and three roll snapshots from each other. Using four clusters, the two-roll snapshots are separated into two clusters; for five clusters, the one-roll snapshots  are additionally separated into two clusters. This information can be gleaned from Fig.~\ref{fig:cavity_3_4_5c}. There, the crosses correspond to full-order solutions color coded to show which k-means snapshot cluster mean it is closest to, i.e., the colors correspond the different k-means clusters and therefore to different local ROM bases. 

\vskip5pt
\noindent\textbf{\em Assigning new parameters to a local basis.} 
As before, we compare two approaches to the assignment problem: the distance to the parameter cluster mean and the distance to mid-range/cluster radius criteria. For both criteria, we then determine the bifurcation diagram using the local ROMs; see Fig.~\ref{fig:cavity_3_4_5c}. As mentioned above, the crosses in that figure correspond to a few full-order solutions; the solid curves are determined from 100 local ROM solutions. 

From the left-column plots in Fig.~\ref{fig:cavity_3_4_5c}, we observe that the parameter cluster mean criterium mis-assigns new parameters in the range Gr$\in$ $[ 80 \cdot 10^{3}, 90 \cdot 10^{3} ]$ to a three-roll basis instead of the appropriate two-roll basis, resulting in a large approximation error. In the three cluster case, the local ROM of one roll solutions is also mistakingly used for two roll solutions corresponding to Gr$\in$ $[ 25 \cdot 10^{3}, 35 \cdot 10^{3} ]$. Other than these mis-assignments the bifurcation diagram is recovered accurately.

From the right-column plots in Fig.~\ref{fig:cavity_3_4_5c}, we observe that the
midrange/radius criterion does not make such mis-assignments and only visibly deviates from the correct bifurcation diagram for the three cluster case in the range
Gr $\in [ 80 \cdot 10^{3}, 90 \cdot 10^{3} ]$. The likely reason is that the problem has not been sampled densely enough in that range.

\vskip5pt
\noindent\textbf{\em Errors and local basis sizes.}
For the cavity, $5$ reference solutions were computed. 
Table~\ref{cavity_7clust_accu} shows the relative errors and basis sizes for $2$, $3$, $5$, and $7$ clusters, giving the same information as the corresponding tables for the jet problem. The local basis sizes are ordered as they are applied with increasing Grashof number. 

\begin{table}[h!]\caption{For the cavity problem: cluster accuracy and basis size comparisons.} 
\begin{tiny}
\begin{tabular}{l*{8}{c}r}
\multicolumn{8}{l}{\underline{\bf $2$ clusters}}\\[.5ex]
ROM   & Gr=12e3 & Gr=32e3 & Gr=52e3 & Gr=75e3 & Gr=97e3 & mean & max & basis size  \\
\hline
Global-1 & 0.0027 & 0.0036 & 0.547 & 0.667 & 0.828   & 0.41 & 0.828 & 18  \\
Global-2 & 0.0029 & 0.0024 & 0.114 & 1.34  & 0.80     & 0.45 & 1.34 & 15  \\
Local    & 0.0028 & 0.0017 & 0.0017 &  0.0015/2.15 & 0.00018 & 0.36 & 2.15  & 15/3  \\[1ex]
\multicolumn{8}{l}{\underline{\bf $3$ clusters}}\\[.5ex]
ROM    & Gr=12e3 & Gr=32e3 & Gr=52e3 & Gr=75e3 & Gr=97e3 & mean & max & basis size  \\
\hline
Global-1 & 0.0027 & 0.0036 & 0.547 & 0.667 & 0.828   & 0.41 & 0.828 & 18  \\
Global-2 & 0.0021 & 0.0035 & 0.0036 & 0.49  & 1.03     & 0.31 & 1.03 & 11  \\
Local    & 0.0028 & 2.177 / 0.0017 & 0.0017 &  0.0015 & 0.00018 & 0.36 & 2.177  & 4/11/3  \\[1ex]
\multicolumn{8}{l}{\underline{\bf $5$ clusters}; $* =$ 2/2/3/6/3}\\[.5ex]
ROM   & Gr=12e3 & Gr=32e3 & Gr=52e3 & Gr=75e3 & Gr=97e3 & mean & max & basis size  \\
\hline
Global-1 & 0.0024 & 0.0027 & 0.722  & 0.983  & 1.086   & 0.56 & 1.086 & 16  \\
Global-2 & 0.131  & 0.0385 & 0.0114 & 0.007  & 0.0004  & 0.038 & 0.131 & 6  \\
Local    & 0.0056 & 0.0055 & 0.5365 & 0.0031 & 0.00018 & 0.11 & 0.5365  &  * \\[1ex]
\multicolumn{8}{l}{\underline{\bf $7$ clusters}; $** =$ 1/3/2/2/2/2/3}\\[.5ex]
ROM   & Gr=12e3 & Gr=32e3 & Gr=52e3 & Gr=75e3 & Gr=97e3 & mean & max & basis size  \\
\hline
Global-1 & 0.0029 & 0.0024 & 0.114  & 1.34  & 0.80   & 0.45 & 1.34 & 15  \\
Global-2 & 0.255  & 2.83  & 0.967 & 2.41  & 0.238  & 1.34 & 2.83 & 3  \\
Local    & 0.0587 & 0.0050 & 0.0090 & 0.013/0.0049 & 0.00018 & 0.015 & 0.0587  & **  \\
\end{tabular}
\end{tiny}
\label{cavity_7clust_accu}
\end{table}

For two clusters, the switching between local ROMs occurs at Gr $ =91.0 \cdot 10^{3}$ according to the midrange criterion and Gr $= 69.3 \cdot 10^{3}$ according to the parameter mean criterion. For three clusters, it occurs at Gr $ = 36.6 \cdot 10^{3}$ and Gr $ = 75.5 \cdot 10^{3}$ according to the parameter mean criterion. According to the midrange criterion, switching occurs at Gr $ = 26.0 \cdot 10^{3}$ and Gr $ = 91.0\cdot 10^{3}$. 
For five clusters, it occurs at 
Gr $ = 16.1 \cdot 10^{3}, 29.8 \cdot 10^{3}, 55.0 \cdot 10^{3}, 82.3 \cdot 10^{3}$
according to the parameter mean criterion.
According to the midrange criterion, switching occurs at 
Gr $ = 16.0 \cdot 10^{3}, 26.0 \cdot 10^{3}, 52.5 \cdot 10^{3}, 91.0 \cdot 10^{3}$.  
For seven clusters, the switching between local ROMs occurs at 
Gr $ = 14.3 \cdot 10^{3}, 26.8 \cdot 10^{3}, 42.5 \cdot 10^{3}, 57.5 \cdot 10^{3},75.0 \cdot 10^{3},89.8 \cdot 10^{3}$  
according to the parameter mean criterion.
According to the midrange criterion, switching occurs at 
Gr $ = 12.5 \cdot 10^{3}, 26.0 \cdot 10^{3}, 42.5 \cdot 10^{3}, 57.5 \cdot 10^{3}, 75.0 \cdot 10^{3}, 91.0 \cdot 10^{3}$.
Because reference value Gr$ = 75.0 \cdot 10^{3}$ corresponds to a switching point, results for both local ROM approximations are shown.

\subsection{Further discussion of SEM errors for both examples}

A mean accuracy of $2\%$ was achieved by local ROMs in 5 out of 8 test cases when using the parameter mean criterion and in 7 out of 8 test cases when using the midrange rule, including the 2 and 3 cluster cases for the second example. A mean accuracy of $1\%$ was achieved by local ROMs in 3 out of 8 test cases, when using the parameter mean criterion and in 5 out of 8 test cases when using the midrange rule, again including the 2 and 3 cluster cases for the second example.

What goes wrong in the failing cases?

First, for the 3 cluster jet case, there is a large error of $2.8\%$ at Re = 37. This value is close to a switching point (at Re = 38)
and the approximation quality tends to degrade at such points as there is a transition period in which one local ROM leaves its parameter range
of validity, whereas the parameter range of validity of the next local ROM has not yet started. Also a large error of $4.5\%$ occurs at Re = 99.
This is the highest Re value considered and all ROMs have difficulties recovering solutions accurately as the Reynolds number approaches 100.
Remedies for this can be to further increase the local ROM size or add more high Re snapshots to the sampling data. 

Second, for the 5 cluster cavity case, there is a large error of $53\%$ at Gr$ = 52.0 \cdot 10^{3}$.
This value is also close to a switching point which is at Gr$ = 52.5 \cdot 10^{3}$ and Gr$ = 55 \cdot 10^{3}$, respectively, and
the same problem of being ``in between'' two clusters holds.

Third, for the 7 cluster cavity case, there is a large error of $5.9\%$ at Gr$ = 12.0 \cdot 10^{3}$.
In this case, the local ROM has only one basis function, which does not provide good approximations.
A remedy for this would be a higher sampling density, such that ROMs having only one basis function are avoided.

By contrast, the best mean accuracy a global ROM achieved is $3.8\%$ for the 5 cluster cavity case. All other global ROMs are above $10\%$ mean accuracy.
This means that no global ROM could accurately recover all solutions in the parameter domain of interest.

For the jet flow, the sample density is higher around the bifurcation point at approximately Re = 33. An increasing
global basis size improves approximation accuracy around the bifurcation point at the cost of accuracy at higher values of Re.
Each global ROM for the jet flow failed to recover solutions for Re higher than 88. 
Because flow behavior changes more significantly at higher Re, a global ROM would need to incorporate that into its projection space.
Local ROMs, on the other hand, do recover the solutions for Re higher than 88. 

In the cavity example, the fixed point algorithm does converge for the global ROMs, but solutions for Grashof numbers higher than Gr$ = 52.0 \cdot 10^{3}$ are 
completely inaccurate. The clustering always separates solutions with 3 rolls from solutions with 1 and 2 rolls, which makes the local ROM approach work.
In the global ROMs, solutions with 1, 2, and 3 rolls are all represented in the same projection space, which could explain the inaccurate approximations.

To conclude, a local ROM approach can work where global ROMs fail. 
It is important to be aware that the approximation quality can degrade in transition regions between two clusters and that, in those regions, the full-order 
model has to be sampled fine enough.

\section{Concluding remarks}\label{sec:con}

In the paper, a k-means clustering of snapshots is used as the starting point for constructing a ROM that uses localized bases to treat PDE problems having bifurcating solutions. A recipe for detecting which local basis to use for any given parameter point not used to determine the snapshots is also given. Careful attention given to account for the differences between bifurcations that cause continuous or discontinuous changes in the solution. At this point, for the sake of simplicity, the new methodology is applied to problems for which the bifurcation diagram is known a priori, that is, one knows which region in the bifurcation diagram a parameter belongs to. Although restrictive, there are may problems for which such knowledge is available. The efficiency of the new approach compared to the use of global ROM bases is demonstrated for simple, one-parameter problems for both types of bifurcations. 

The results obtained so far are encouraging so that current and future work focuses on further developing our approach so that it can become a useful tool for practitioners. In particular, we will focus on:

\begin{itemize}

\item[--] testing the methodology for problems with several input parameters;

\item[--] automatically detecting, without having any a priori knowledge, which local basis one should use to a given parameter point; 

\item[--] testing the methodology for more ``realistic'' problems;

\item[--] application of the methodology to control and optimization and to uncertainty quantification problems.

\end{itemize}

\bibliographystyle{plain}
\bibliography{rbsissa,latexbi}

\end{document}